\documentclass[10pt]{article}
\usepackage{etex}
\usepackage[textwidth=15cm, textheight=22cm]{geometry}
\usepackage{amsfonts}
\usepackage{amsmath}
\usepackage{amssymb}
\usepackage{amsthm}
\usepackage{enumitem}
\usepackage{graphicx}

\usepackage[dvipsnames]{xcolor}

\usepackage{tikz}
\usetikzlibrary{arrows.meta,shapes,calc}
\usetikzlibrary{positioning,shapes,arrows}

\usepackage{adjustbox}

\usepackage{pgfplots}
\pgfplotsset{compat=1.12}
\usetikzlibrary{plotmarks}
\usepackage[ngerman, english]{babel}
\usepackage{babelbib}
\usepackage[latin1]{inputenc}

\usepackage{nicefrac}

\usepackage{subfig}
\usepackage{epstopdf}

\usepackage{bbm}

\usepackage{url}

\usepackage[colorinlistoftodos,prependcaption,textsize=tiny]{todonotes}

\usepackage{kpfonts}

\usepackage{algorithmic}
\usepackage[section,Algorithm]{algorithm}

\usepackage{epstopdf}
\usepackage[dvips]{epsfig}

\usepackage{graphicx}%

\newcommand{\AAA}{\mathcal{A}}

\newcommand{\EE}{\mathcal{E}}
\newcommand{\FF}{\mathcal{F}}
\newcommand{\GG}{\mathcal{G}}

\newcommand{\NN}{\mathcal{N}}

\newcommand{\TT}{\mathcal{T}}
\newcommand{\UU}{\mathcal{U}}
\newcommand{\VV}{\mathcal{V}}

\newcommand{\E}{\mathbb{E}}

\newcommand{\R}{\mathbb{R}}

\newcommand{\tg}{\tilde{g}}

\newcommand{\f}{f} 
\newcommand{\n}{2} 
\newcommand{\aae}{a^{e}} 
\newcommand{\be}{b^{e}} 
\newcommand{\rhoe}{q^{e}} 
\newcommand{\fe}{f^{e}} 
\newcommand{\s}{s} 
\newcommand{\B}{B} 
\newcommand{\Ga}{\Gamma_a} 
\newcommand{\Gae}{\Gamma_a^{e}} 
\newcommand{\Gb}{\Gamma_b} 
\newcommand{\Gbe}{\Gamma_b^{e}} 
\newcommand{\ga}{\gamma_a} 
\newcommand{\gae}{\gamma_a^{e}} 
\newcommand{\gb}{\gamma_b} 
\newcommand{\gbe}{\gamma_b^{e}} 
\newcommand{\Gabe}{\Gamma_{a/b}^{e}} 
\newcommand{\gabe}{\gamma_{a/b}^{e}} 
\newcommand{\couplingv}{c_v} 
\newcommand{\pe}{p^{e}} 
\newcommand{\qe}{q^{e}} 
\newcommand{\Ue}{U^{e}} 
\newcommand{\Ie}{I^{e}} 

\newcommand{\EEv}{\EE^{v}} 
\newcommand{\ev}{e^{v}} 
\newcommand{\Uev}{U^{\ev}} 

\newcommand{\Supp}{S_{u(t)}} 
\newcommand{\SuppO}{S_{u(t_0)}} 
\newcommand{\detCost}{C_{u(t)}} 
\newcommand{\compr}{u_{compr}(t)} 

\newcommand{\risk}{\theta} 
\newcommand{\cdfJCC}{G} 

\newcommand{\densGauss}{\varphi} 
\newcommand{\densTGauss}{\rho_{X,a,b}} 
\newcommand{\mTGauss}{\xi} 
\newcommand{\vTGauss}{\sigma^2} 
\newcommand{\sdevTGauss}{\sigma} 

\newcommand{\wdet}{w_1} 
\newcommand{\wtrack}{w_2} 
\newcommand{\wunder}{w_3} 
\newcommand{\wex}{w_4} 

\newcommand{\spaceArc}{[a^{e}, b^{e}]} 

\newcommand{\invLa}{\nicefrac{1}{\lambda}}

\newcommand{\ind}{\mathbbmss{1}}

\newlength{\fwidth}

\DeclareMathOperator{\cov}{Cov}
\DeclareMathOperator{\id}{id}



\theoremstyle{plain}
\newtheorem{thm}{Theorem}[section]
\newtheorem{lem}[thm]{Lemma}
\newtheorem{prop}[thm]{Proposition}

\theoremstyle{definition}
\newtheorem{defin}{Definition}[section]

\theoremstyle{remark}
\newtheorem*{rem}{Remark}

\begin{document}
\title{Chance-constrained optimal inflow control in hyperbolic supply systems with uncertain demand}
\author{Simone G\"{o}ttlich\footnotemark[1], \; Oliver Kolb\footnotemark[1], \; Kerstin Lux\footnotemark[1]}
\footnotetext[1]{University of Mannheim, Department of Mathematics, 68131 Mannheim, Germany (goettlich@uni-mannheim.de, klux@mail.uni-mannheim.de).}
\maketitle

\begin{abstract}
In this paper, we address the task of setting up an optimal production plan ta\-king into account an uncertain demand. 
The energy system is represented by a system of hyperbolic partial differential equations (PDEs) and the uncertain demand stream is captured by an Ornstein-Uhlenbeck process. We determine the optimal inflow depending on the producer's risk preferences. The resulting output is intended to optimally match the stochastic demand for the given risk criteria. We use uncertainty quantification for an adaptation to different levels of risk aversion. More precisely, we use two types of chance constraints to formulate the requirement of demand satisfaction at a prescribed probability level.
In a numerical analysis, we analyze the chance-constrained optimization problem for the Telegrapher's equation and a real-world coupled gas-to-power network.
\end{abstract}

\section{Introduction} \label{sec:Intro}

In recent years, significant attention has been paid to the energy market. On the one hand, this is due to climate protection policies. On the other hand, the decision of the German government on the nuclear phase-out by the end of 2022\footnote{\url{https://www.bmwi.de/Redaktion/EN/Artikel/Energy/nuclear-energy-nuclear-phase-out.html
}, last checked: $4^{th}$ of April, 2019} will cause significant changes in the energy sector. Those changes are to a large extent triggered by the specification to have $65\%$ of the gross electricity consumption generated out of renewable energy sources by the end of 2030. One major challenge is the handling of uncertainty in the renewable energy production. It heavily depends on the weather conditions and is subject to large fluctuations. Those fluctuations can heavily affect the grid operation or even lead to outages (\cite{Bienstock.2014}).

Another source of uncertainty in the power grid operation are power demands (see \cite{Zhang.2011}).
This makes it difficult to guarantee a sufficient supply thereby influencing reliability and profitability in power system operation (see \cite{Geletu.2013}). It is crucial for all players in the electricity sector to cope with stochastic demand fluctuations.
In the context of gas-to-power, those fluctuations might carry over to the gas network being coupled to the electricity grid. The amount of gas converted to power and feeded into the electricity system depends on the power demand and therefore inherits its stochasticity. Furthermore, the gas demand itself is uncertain (see e.g.\ \cite{Henrion.2015}).

\bigskip

There are several ways to formulate optimization problems in the presence of uncertainty. As in \cite{Henrion.2017}, we distinguish between \textit{robust optimization} and \textit{chance constrained optimization (CCOPT)}. A robust formulation is appropriate if the distribution of the uncertain parameter in the system is unknown because it cannot be measured or observed (no historical data), whereas a probabilistic programming (CCOPT) is suitable in the presence of historical observations where a distribution of the unknown quantity can be derived (see \cite{Henrion.2017}).

In this paper, we focus on uncertain demands. In this context, it is reasonable to assume access to historical demand data and we therefore focus on CCOPT for the remainder of the present manuscript. Chance constraints have been introduced in 1958 in \cite{Charnes.1958} by Charnes, Cooper and Symonds in the context of production planning, and further developed in \cite{Charnes.1959}. A lot of work has been done in CCOPT. It is worthwhile to particular mention the contributions of Pr\'{e}kopa. His monograph on stochastic programming from 1995 \cite{Prekopa.1995} nowadays still serves as a standard reference. A general overview on properties, solution methods, and fields of application such as hydro reservoir management (\cite{Henrion.2010,Henrion.2014}), optimal power flow (\cite{Bienstock.2014,Zhang.2011}), energy management (\cite{vanAckooij.2013}), and portfolio optimization (\cite{Rockafellar.2000}) can be found in \cite{Geletu.2013}.

\bigskip
In this work, we are interested in a constrained optimization problem composed of three constraints: a stochastic differential equation (SDE) that models the uncertain demand, a system of hyperbolic balance laws to describe the energy system (electricity or gas), and a chance constraint ensuring demand satisfaction with a certain probability (risk level).
Chance constraints in the context of PDE-constrained optimization are also considered in \cite{Surowiec.2018}. However, the requirement of a continuously Fr\'{e}chet-differentiable solution of the PDE rules out most hyperbolic PDEs due to the possibility of discontinuous solutions even for continuously Fr\'{e}chet-differentiable initial data.

In \cite{Henrion.2018}, they apply their results on semi-continuity, convexity, and stability in an infinite-dimen\-sional setting to a simple PDE-constrained control problem. Again, the PDE is not of hyperbolic nature.

\bigskip
There already exist a few investigations of CCOPT for systems of hyperbolic nature such as gas networks. Accounting for the stochasticity of demand, \cite{Henrion.2016} presents a method to compute the probability of feasible loads in a stationary, passive gas network. Note that the steady state assumption on the network entails constraints formulated in terms of purely algebraic equations instead of full hyperbolic dynamics. In \cite{Henrion.2017}, they extend the aforementioned work by also allowing for uncertainty with respect to the roughness coefficient. They analyze the question of the maximal uncertainty allowed such that random loads can be satisfied at a prescribed probability level within a stationary, passive gas network. The feasibility of random loads, again in a stationary gas model, this time with active elements in terms of compressor stations is considered in \cite{Gugat.2018}. In \cite{Adelhuette.2017}, they also account for the transient case. However, they do not take into account the full hyperbolic dynamic but the linear wave equation instead.

In contrast to the above-mentioned contributions, here, we consider the full isentropic gas dynamics. Note however that the full stochasticity of the demand in the power-to-gas setting presented in Section \ref{subsec:OCgasPower} does not enter the nonlinear dynamics of the gas flow directly but is coupled to it via the objective function.

Another important difference to the aforementioned contributions is the modelling of the uncertain demand. In many cases, (truncated) Gaussian random vectors are used as proposed in \cite{Henrion.2015} and used in for example \cite{Adelhuette.2017,Henrion.2017,Gugat.2018}. Using stochastic processes to model the demand enables to capture time dynamics. Sometimes the demand is modeled by a discrete time stochastic process as e.g.\ in \cite{Henrion.2010}. In \cite{Henrion.2014}, they use a continuous time stochastic process additively decoupled in a deterministic trend and a causal process generated by Gaussian innovations. In contrast to that, we consider a dynamic demand model for the random loads via a continuous time stochastic process, the so-called Ornstein-Uhlenbeck process, where the deterministic trend and the stochastic evolution of the process are coupled. This is also a popular choice to model demand in various fields (see e.g.\ \cite{Barlow.2002} in the context of electricity). In \cite{Nesti.2016}, and \cite{Wadman.2016}, they use a multivariate Ornstein-Uhlenbeck process from the supply-side point of view. They model uncertain injections into the power network and assess the probability of outages. A desirable feature for both supply and demand modeling is the mean reverting property. The process is always attracted to a certain predefined mean level, which may result from historical supply/demand data. In contrast to those contributions, we use a time-dependent mean level, which enables us to depict seasonal and daily patterns appearing in historical data also in our mean-demand level (no longer assumed to be constant) and set up an optimization framework based thereon.

\bigskip
The paper is organized as follows: in Section \ref{sec:SOC}, we introduce the stochastic optimal control setting. We define the energy system as well as the stochastic demand process, and introduce the considered cost functional. We distinguish between a single chance constraint (SCC) as well as a joint chance constraint (JCC). In Section \ref{sec:reformStochComp}, we consider a deterministic reformulation of both the cost functional and the SCC, and present a stochastic reformulation of the JCC, which allows to incorporate it in our numerical framework. In Section \ref{sec:numRes}, we validate the numerical framework for a simple case of a hyperbolic supply system, i.e. the scalar linear advection with source term. It represents a suitable test case setting as we are able to derive the corresponding analytical solution in case of an SCC. We then apply our numerical routine to the Telegrapher's equations, a linear system of hyperbolic balance laws. We conclude with a numerical investigation of a nonlinear system of hyperbolic balance laws in terms of a real-world gas-to-power application.

\section{Stochastic optimal control setting} \label{sec:SOC}
In this section, we set up the mathematical framework for the task of finding an optimal injection plan taking into account an uncertain demand for the time period from $t_0=0$ to final time $T$. This has been done for a linear transport equation in \cite{Lux.2018}. Here, we use the stochastic optimal control framework set up in \cite{Lux.2018}, and extend it to more complex supply dynamics on a network. We model a general energy system by a system of hyperbolic balance laws on a network, and the stochastic demand is described by an Ornstein-Uhlenbeck process (OUP).

\bigskip

\subsection{Energy system with uncertain demand} \label{subsec:hypSupplySystem}

We consider different types of $\n$-dimensional energy systems. The network is modeled by a finite, connected, directed graph $\GG=(\VV,\EE)$ with a non-empty vertex (node) set $\VV$ and a non-empty set of edges $\EE$. For $v \in \VV$, we define the set of all incoming edges by $\delta^{-}(v) = \{e \in \EE: e=(\cdot,v) \}$, and the set of all outgoing edges by $\delta^{+}(v)=\{e \in \EE: e=(v,\cdot) \}$.
In the sequel, we identify an edge $e=(v_{in},v_{out})$ by the interval $[a_e,b_e]$, where $a_e$ denotes the starting point of the edge, and $b_e$ its end point.
We further define the set of inflow vertices by $\VV_{in}=\{v\in \VV: \delta^{-}(v) = \emptyset \}$, and the set of outflow vertices by $\VV_{out}=\{v\in \VV: \delta^{+}(v) = \emptyset \}$. As a simplification, here, we restrict our network to $|\VV_{in}|=|\VV_{out}|=1$. Note however that an extension to several inflow nodes ($|\VV_{in}|>1$), and outflow nodes ($|\VV_{out}|>1$) is straightforward by considering a vector-valued inflow control, and a multivariate Ornstein-Uhlenbeck process as used in \cite{Wadman.2016}.

The inflow control $u(t)$ acts on the vertex $v_{in} \in \VV_{in}$, and the demand $Y_t$ is realized at the vertex $v_{d} \in \VV_{out}$. We require $b^{e}=L$ for all $e \in \delta^{-}(v_{d})$.

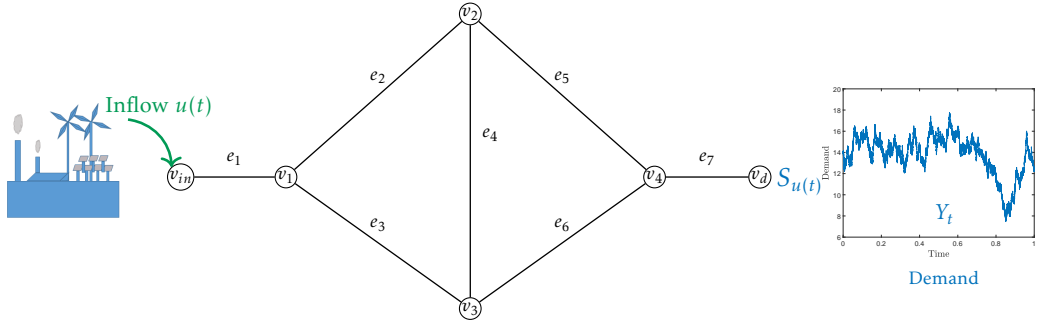
\begin{figure}[h!]
	\centering
	\setlength{\fwidth}{0.3\textwidth}
	\begin{tikzpicture}[transform shape,scale=0.7]
	\node (powerStation) {\includegraphics[width=0.15\linewidth]{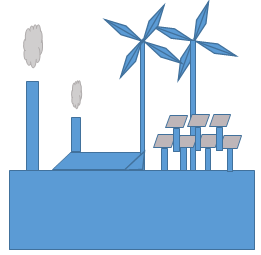}};
	\node[circle, draw=black, fill=white, inner sep=0pt] (xIn) at (powerStation.east)[xshift =1cm ,yshift=-0.3cm]{$v_{in}$};
	\node[circle, draw=black, inner sep=0pt] (x0) at (powerStation.east)[xshift =3cm ,yshift=-0.3cm]{$v_1$};
	\node[circle, draw=black, inner sep=0pt] (xEnd) at (powerStation.east)[xshift =10cm ,yshift=-0.3cm]{$v_4$};
	\node[circle, draw=black, inner sep=0pt] (xUp) at (powerStation.east)[xshift =6.5cm ,yshift=2.8cm]{$v_2$};
	\node[circle, draw=black, inner sep=0pt] (xDown) at (powerStation.east)[xshift =6.5cm ,yshift=-2.8cm]{$v_3$};
	\node[circle, draw=black, fill=white, inner sep=0pt] (xDem) at (powerStation.east)[xshift =12cm ,yshift=-0.3cm]{$v_{d}$};
	\node at (xDem.north)[xshift =0.75cm ,yshift=-0.3cm]{\Large \textcolor{RoyalBlue}{$\Supp$}};
	\node (demand)  at (powerStation.east)[xshift =15.2cm ,yshift=-0.2cm]{\includegraphics[width=0.27\linewidth]{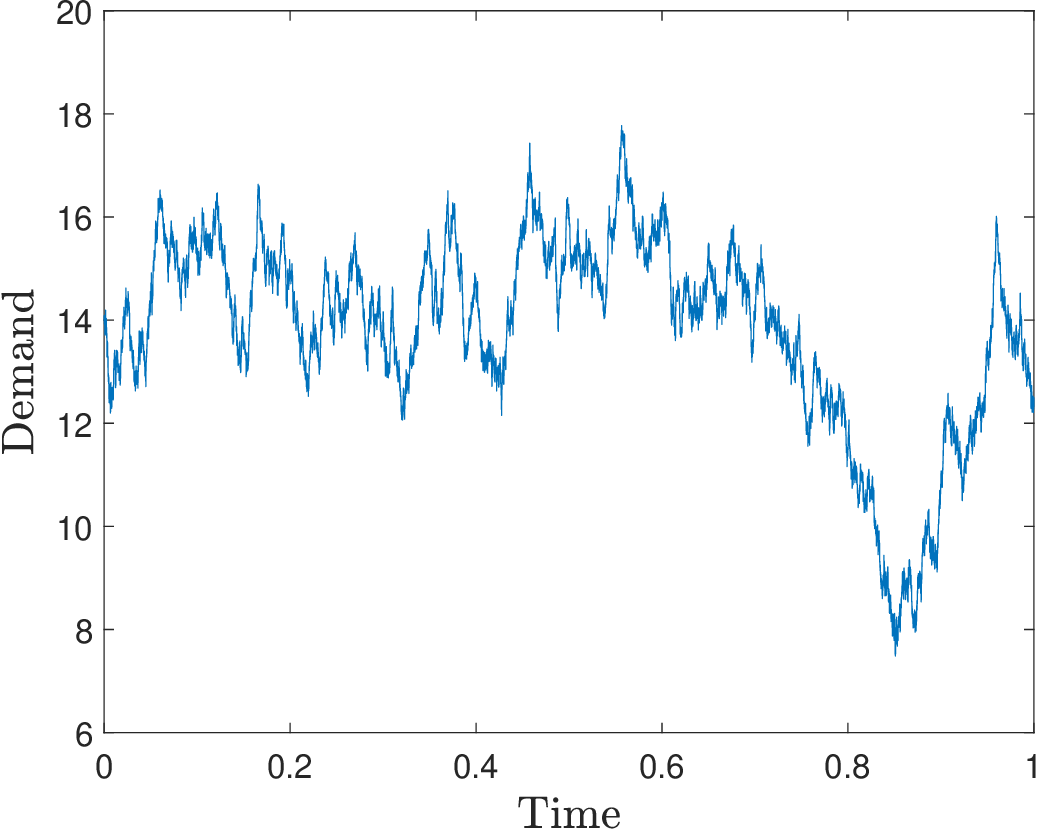}};
	\node at (powerStation.east)[xshift =15.5cm ,yshift=-2.2cm]{\textcolor{RoyalBlue}{Demand}};
	\node (demand)  at (powerStation.east)[xshift =15.5cm ,yshift=-1cm]{\large \textcolor{RoyalBlue}{$Y_t$}};
	\draw [->, line width=1pt, ForestGreen] (powerStation) to[bend left=40]
	node[above,yshift=0.2cm] {\large Inflow $u(t)$} (xIn);
	\draw [-, line width=0.5pt, black] (xIn) to node[above,yshift=0.1cm] {$e_1$} (x0);
	\draw [-, line width=0.5pt, black] (x0) to node[above,yshift=0.1cm] {$e_2$}	(xUp);
	\draw [-, line width=0.5pt, black] (xUp) to node[above,yshift=0.1cm] {$e_5$} (xEnd);
	\draw [-, line width=0.5pt, black] (x0) to node[above,yshift=0.1cm] {$e_3$}	(xDown);
	\draw [-, line width=0.5pt, black] (xDown) to node[above,yshift=0.1cm] {$e_6$}	(xEnd);
	\draw [-, line width=0.5pt, black] (xUp) to node[right,xshift=0.1cm, yshift=0.5cm] {$e_4$} (xDown);
	\draw [-, line width=0.5pt, black] (xEnd) to node[above,yshift=0.1cm] {$e_7$} (xDem);
	\end{tikzpicture}
	\caption{Example of an energy network with uncertain demand at $v_d$}
	\label{fig:networkTop}
\end{figure}

The dynamics of the energy system on one edge $e$ are given by the following hyperbolic balance law with initial condition (IC) and boundary conditions (BCs)
\begin{align}
\partial_t \rhoe + \partial_x \fe(\rhoe) &= \s(\rhoe), \quad x \in \spaceArc, \ t \in [0,T] \notag\\
\rhoe(x,0) &= q_0^{e}(x), \notag\\
\Gae(\rhoe(a_e,t)) &= \gae(t), \quad \Gbe(\rhoe(b_e,t)) = \gbe(t). \label{eq:IBVP} 
\end{align}
Thereby, $\f:\R^{\n} \rightarrow \R^\n$ is a given flux function and $\s:\R^{\n} \rightarrow \R^{\n}$ the source term.
$\rho_0^{e}: \spaceArc \rightarrow \R^{\n}$ describes the initial state of the system on edge $e$. The functions $\Gabe:\R^\n \rightarrow \R^{m_{l/r}}$ enable to prescribe a certain evaluation of a certain number of components of the density at the boundary. Thereby, $0 \leq m_l,m_r \leq \n$ denote the number of prescribed BCs at the left respectively right boundary. A value of $m_{l/r}=0$ has to be interpreted in a way that no left/right boundary condition is prescribed. For a scalar conservation law, choosing $m_l=1$, $m_r=0$, and $\Gae=\f$ corresponds to prescribing only the flow at the left boundary.

The boundary conditions (BCs) themselves are given in terms of the functions $\gae: [0,T] \rightarrow \R^{m_l}$ at the left boundary $a_e$ of edge $e$, and $\gbe: [0,T] \rightarrow \R^{m_r}$ at the right boundary $b_e$ of edge $e$.
The numbers of BCs $m_{l/r}$, and the functions $\gabe$ need to be chosen carefully such that they are consistent with the characteristics of the conservation law. This will be further specified below for each setting \ref{subsec:OCtelEq}-\ref{subsec:OCgasPower}.

Moreover, for the non-degenerate network case, i.e. $\exists v \in \VV: |\delta^{-}(v)|+|\delta^{+}(v)|\geq 2$, suitable coupling conditions $\couplingv: \R^{\n \times \left(|\delta^{+}(v)|+|\delta^{-}(v)|\right)} \rightarrow \R^{m_c}$ for each vertex need to be imposed. Thereby, $m_c$ denotes the number of coupling conditions. This will be made explicit in Subsection \ref{subsec:OCgasPower}.

For $e \in \delta^{+}(v_{in})$, the function $\gae(t)$ depends on the inflow control $u(t)$.

To simplify notation, we do not explicitly write down the dependence of the supply on the density at the end points of all ingoing edges and denote the supply at $v_d$ at time $t$ simply by $\Supp$ meaning
\begin{align*}
	\Supp=S_u\left(\varprod_{e \in \delta^{-}(v_d)} \rhoe(b_e,t)\right).
\end{align*}


\subsubsection*{Discretization scheme}
For the numerical investigation in Section \ref{sec:numRes}, the considered hyperbolic energy systems need an appropriate discretization scheme. Motivated by our real-world example, we choose an implicit box scheme (IBOX)~\cite{KolbLangBales2010} for all considered scenarios. For a general system of balance laws (on any edge)
\begin{align*}
\partial_t q + \partial_x f(q) = s(q),
\end{align*}
the considered scheme reads
\begin{align}
\frac{Q^{n+1}_{j-1} + Q^{n+1}_{j}}{2} = \frac{Q^{n}_{j-1} + Q^{n}_{j}}{2} - \frac{\Delta t}{\Delta x}\left( f(Q^{n+1}_j)- f(Q^{n+1}_{j-1}) \right) + \Delta t \frac{s(Q^{n+1}_j)+s(Q^{n+1}_{j-1})}{2}. \label{eq:IBOXscheme}
\end{align}
Here, $\Delta t$ and $\Delta x$ are the temporal and spatial mesh size, respectively, and the numerical approximation is thought in the following sense:
\begin{align}
Q_j^n \approx q(x,t) \quad \text{for} \quad x\in X_j=\big[ a+(j-\tfrac{1}{2})\Delta x , a+(j+\tfrac{1}{2})\Delta x \big) \cap \big[ a , b \big] , \ t \in I_i=\big[ i\Delta t , (i+1)\Delta t \big).
\end{align}
%
To avoid undesired boundary effects, the discretization of the initial condition on bounded domains is done pointwise, i.e.,
\begin{align}
 Q_j^n = q_0(a + j \Delta x).
 \label{eq:IBOX_IC}
\end{align}
As remarked in \cite{KolbLangBales2010}, for a discretization $x_l<x_{l+1}<\cdots<x_{r-1}<x_r$, we obtain $r-l$ equations for $r-l+1$ variables (in the scalar case). This entails the need to prescribe BCs at exactly one boundary specified by the characteristic direction. This also explains the above mentioned assumption of no change in the signature of the characteristic directions on the considered domain.
The discrete version of the BCs is given by
\begin{align}
\Ga(Q_0^i) = \ga(t_i), \quad \Gb(Q_{N_{\Delta x}}^{i}) = \gb(t_i) \quad \forall \ i \in \{1,\cdots,N_{\Delta t}\}. \label{eq:IBOX_BCs}
\end{align}

\bigskip
Note that the implicit box scheme has to obey an inverse CFL condition~\cite{KolbLangBales2010}, which is beneficial for problems with large characteristic speeds whereas the solution is merely quasi-stationary. This is usually the case for daily operation tasks in gas networks and therefore motivates the choice within this work.


\subsubsection*{Uncertain Demand} \label{subsubsec:uncertainDemand}
The uncertain demand stream is modeled by an Ornstein-Uhlenbeck process (OUP), which is the unique strong solution of the following stochastic differential equation
\begin{align}
dY_t = \kappa\left(\mu(t)-Y_t\right)dt + \sigma dW_t,\ \quad Y_{t_0}=y_{t_0} \label{eq:OUP}
\end{align}
on the probability space $\left(\Omega,\FF,P\right)$. $W_t$ is a given one-dimensional Brownian motion on the same probability space, and $y_0$ describes the demand at time $t_0=0$. The constants $\sigma>0$ and $\kappa>0$ describe the speed of mean reversion and the intensity of demand fluctuations. By mean reversion, we refer to the property of the process that it is always attracted by a certain time-dependent level $\mu(t)$, called the mean demand level. This is due to the sign of the drift term $\kappa\left(\mu(t)-Y_t\right)$, which ensures that being above (below) the mean demand level, the process experiences a reversion back to it.
By the OUP \eqref{eq:OUP}, we are able to capture daily or seasonal patterns in the demand.
From a mathematical point of view, this process has some nice analytical properties. For example, we can derive the solution to equation \eqref{eq:OUP} explicitly via the It\^{o} formula. It reads as
\begin{align}
Y_t &= y_{t_0} e^{-\kappa (t-t_0)} + \kappa \int_{t_0}^t{e^{-\kappa(t-s)} \mu(s)ds} + \sigma \int_{t_0}^t{e^{-\kappa(t-s)}dW_s}. \label{eq:expSolOUP}
\end{align}
Moreover, it is possible to derive its distribution explicitly as
\begin{align}
Y_t &\sim N {\underbrace{\left(y_{t_0}e^{ - \kappa (t-t_0)}  + \kappa \int\limits_{t_0}^t {e^{ - \kappa \left( {t - s} \right)} \mu \left( s \right)ds}\right.}_{m_{OUP}(t)} ,\,\,\underbrace{\left. \vphantom{y_{t_0}e^{ - \kappa (t-t_0)}  + \kappa \int\limits_{t_0}^t {e^{ - \kappa \left( {t - s} \right)} \mu \left( s \right)ds}}\frac{\sigma ^2}{2\kappa}\left(1 - e^{ - 2\kappa \left( {t - t_0} \right)}\right)\right)}_{v_{OUP}(t)} \,}. \label{eq:distOUP}
\end{align}
For further details on the demand process, we refer the reader to \cite{Lux.2018}.

\subsection{Chance constraints} \label{subsec:CCs}
Requiring demand satisfaction for every realization of the demand process might be too restrictive as, in some cases, it might lead to an infeasible optimization problem (\cite[p.\ 5]{Shapiro.2009}).

One possibility to overcome problems of infeasibility and to reduce the average undersupply is to introduce an undersupply penalty term in the cost function. The effect of an undersupply penalty on the optimal output has been analyzed in \cite{Lux.2019}. A comparison of different types of undersupply can be found in \cite{Lux.2019b}.

\bigskip 
Another approach is to guarantee with a certain probability that there is no undersupply within a prescribed time interval $I_{CC} \subset [t^{*},T]$, where $t^{*}$ is the first time that a supply is realized at $v_d$. Mathematically this is formulated in terms of a chance constraint (CC). One possibilty is to require \textit{at each point in time} that the probability of a demand satisfaction is at least equal to one minus a given risk level $\risk$ (see \eqref{eq:SCC}). This results in a so called single chance constraint (SCC). Another possibility is a joint chance constraint (JCC), that is, we require that the probability of a demand satisfaction is at least equal to one minus a given risk level $\risk$ \textit{on a whole interval} simultaneously (see \eqref{eq:JCC}).
\begin{subequations} \label{eq:CC}
	\begin{align}
	P\left(Y_t \leq \Supp\right) &\geq 1-\risk \ \forall \ t \in I_{CC}, \label{eq:SCC} \\
	P\left(Y_t \leq \Supp \ \forall \ t \in I_{CC}\right) &\geq 1-\risk. \label{eq:JCC} 
	\end{align}
\end{subequations}

\subsection{Objective function and stochastic optimal control problem} \label{subsec:OF}
Having formulated all the optimization constraints, we now address the objective function.
We aim at minimizing the arising costs. We construct our cost function out of several components. One component consists of deterministic costs $\detCost$ such as operating costs for a gas compressor ($C1$).
Another component are tracking type costs that arise from a mismatch between the externally given demand \eqref{eq:OUP} and our supply $S_u$ realized at the demand vertex $v_d$. We measure the tracking type costs in terms of the expected quadratic deviation between the demand and supply ($C2$).
Particularly accounting for negative mismatches, we introduce undersupply costs as a third cost component ($C3$).
If a sale of excess supply is possible, the revenue can be included into the cost function as a fourth component ($R$).

Putting the deterministic costs $C1$, the tracking costs $C2$, the undersupply costs $C3$, and the excess revenue $C4$ together, we obtain a cost function of the following type:
\begin{align}
OF(Y_t, t_0,y_{t_0} ,\Supp) =& \underbrace{\vphantom{\E\left[\left(\Supp-Y_t\right)^2|Y_{t_0}=y_{t_0}\right]} \wdet \cdot \detCost}_{C1} + \underbrace{\wtrack \cdot \E\left[\left(\Supp-Y_t\right)^2|Y_{t_0}=y_{t_0}\right]}_{C2} \notag\\
&- \underbrace{\vphantom{\E\left[\left(\Supp-Y_t\right)^2|Y_{t_0}=y_{t_0}\right]} \wunder \cdot \E\left[\left(\Supp - Y_t\right)_{-}|Y_{t_0}=y_{t_0}\right]}_{C3} - \underbrace{\vphantom{\E\left[\left(\Supp-Y_t\right)^2|Y_{t_0}=y_{t_0}\right]} \wex \cdot \E\left[\left(\Supp - Y_t\right)_{+}|Y_{t_0}=y_{t_0}\right]}_{R} , \label{eq:generalOF}
\end{align}
where
\begin{align*}
	\left(\Supp - Y_t\right)_{-} &= \left\{\begin{array}{ll}
	\Supp - Y_t & \text{if} \ \Supp<Y_t \\
	0 & \text{if} \ \Supp\geq Y_t
	\end{array}
	\right.
\end{align*}

\bigskip

All together, we come up with the stochastic optimal control (SOC) problem
\begin{align}
\min_{u \in \UU_{ad}} \int_{t^{\ast}}^{T} OF(Y_t; t_0;y_{t_0} ;\Supp) dt \quad \text{subject to} \ \eqref{eq:IBVP}, \eqref{eq:OUP},\ \text{and} \ \eqref{eq:CC}. \label{eq:SOC-CC}
\end{align}

To determine the optimal control, we need to specify measurability assumptions on the control by defining the space of admissible controls $\UU_{ad}$.
\begin{align}
	\UU_{ad} &= \{u:[t_0,T] \rightarrow \R \ | \ u \in L^2\left([t_0,T]\right), u(t) \geq 0, \ \text{and} \ u(t) \ \text{is} \ \FF_{t_0}\text{--predictable} \ \text{for} \ t \in [t_0,T]\}. \label{eq:controlSpaceCM1}
\end{align}
Results for this control method and two other control methods with an objective function of pure tracking type ($\wdet=\wunder=\wex=0$, $\wtrack=1$) for the linear advection equation without imposing CCs can be found in \cite{Lux.2018}.

\section{Deterministic reformulation of the stochastic problem} \label{sec:reformStochComp}

Having set up the SOC problem \eqref{eq:SOC-CC}, the question of how to solve this minimization problem naturally arises. One way is to trace back the SOC problem \eqref{eq:SOC-CC} to a deterministic setting, in which we can apply well-known methods from deterministic PDE-constrained optimization such as adjoint calculus. In order to do so, in \ref{subsec:reformCC}, we analytcially treat both types of CCs presented in Subsection \ref{subsec:CCs}. The deterministic expression of the objective function introduced in \ref{subsec:OF} is derived in \ref{subsec:reformOF}.

\subsection{Reformulation of chance constraints} \label{subsec:reformCC}
The reformulation of the CC heavily depends on the type of CC. Whereas the SCC \eqref{eq:SCC} can be reformulated via quantiles of a normal distribution, the reformulation of the JCC \eqref{eq:JCC} is not obvious. Therefore, we need to treat the different types of CCs separately.

\textbf{Single chance constraint}

For the SCCs \eqref{eq:SCC}, we use a quantile-based reformulation as mentioned in \cite[p.\ 580]{Henrion.2010}. By using the known distribution \eqref{eq:distOUP} of the OUP, this results in the deterministic state constraints
\begin{align}
\Supp \geq m_{OUP}(t) + v_{OUP}(t)\Phi^{-1}\left(1-\risk\right) \ \forall t \in I_{CC}, \label{eq:stateConstr}
\end{align}
where $m_{OUP}(t)=y_{t_0}e^{-\kappa (t-t_0)}  + \kappa \int\limits_{t_0}^t {e^{-\kappa \left( {t - s} \right)} \mu \left( s \right)ds}$, $v_{OUP}(t)=\sigma ^2 \int\limits_{t_0}^t {e^{ - 2\kappa \left( {t - s} \right)}ds}$, and $\Phi$ is the standard normal cumulative distribution function.

\textbf{Joint chance constraint}

JCCs \eqref{eq:JCC} are mathematically by far more involved than SCCs (see \cite[p.\ 1213]{Geletu.2013}). No longer considering the constraint pointwise in time, we now have to deal with a joint probability distribution. Unfortunately, we can no longer make use of the deterministic reformulation as state constraint as for \eqref{eq:SCC}.

As the integration of the JCC into the optimization framework is a core issue to tackle the SOC problem \eqref{eq:SOC-CC}, we need to come up with a different approach.
As in \cite{Wadman.2016}, we now use a non-deterministic reformulation of the JCC as a first passage time problem. We denote by
\begin{align*}
	\TT_{y_0} &= \inf_{t\geq t_0}\left(t \ |\ Y_t > \Supp\right)
\end{align*}
for $Y_{t_0}=y_{t_0}<\SuppO$ the first passage time of the OUP, and obtain the equivalent formulation of the JCC \eqref{eq:JCC} as a first passage time problem of the form
\begin{align}
	P(\TT_{y_0} > t) &\geq 1-\risk. \label{eq:reformJCCasFPT}
\end{align}
For further details on first passage times, we refer the reader to \cite[Chapter 4]{Dynkin.1965}. We will see that reformulation \eqref{eq:reformJCCasFPT} enables to include constraint \eqref{eq:JCC} into our SOC framework.

\bigskip
However, this is not obvious. The drawback is that we have to deal with the distribution of the first passage time of the OUP with time-dependent mean demand level for a time-dependent absorbing boundary. The time-dependency rules out some classical approaches. Furthermore, even for a constant boundary, the task turns out to be by far more complicated than deriving the distribution of the first passage time of a Brownian motion (see \cite{Lipton.2018,Uhlenbeck.1945}).
A closed-form solution is only known for a few special cases.

One reason making the present case particularly hard is that there is no pure diffusion equation any more that would allow to apply the formula presented in the paper \cite{Darling.1953}. So far, to the best of our knowledge, no closed-form solution of the first passage time density in our case is known. Several semi-analytical approaches exist: In \cite{Jiang.2019}, they derive an integral representation of the first-passage time of an inhomogeneous OUP with an arbitrary continuous time-dependent barrier and extend their results to continuous Markov processes. The idea of exploiting transformations among Gauss-Markov processes to relate the problem to the known first passage time density of a Wiener process. However, as stated in \cite{DiNardo.2001}, the transformation of the OUP to the Brownian motion entails exponentially large times. In an alternative approach presented in \cite{DiNardo.2001}, they consider this broader class of real continuous Gauss-Markov process with continuous mean and covariance functions. The latter approach is the one that we tailor to our time-dependent OUP.

In a first step, we show that the approach in \cite{DiNardo.2001} is applicable to our case of the time-dependent OUP, and in a second step, we introduce the method itself.

\bigskip

\textbf{Verification of prerequisites}

As we want to apply a result on the first passage time density formulated in a general Gauss-Markov setting, we first need to verify that the time-dependent OUP given by equation \eqref{eq:OUP} is indeed a Gauss-Markov process.
\begin{lem}
	The OUP given by equation \eqref{eq:OUP} is a Gauss-Markov process.
\end{lem}
\begin{proof}
	We first verify that the OUP is a Gauss process, i.e.\ that for any integer $n\geq 1$ and times $0\leq t_1<t_2<\cdots<t_n \leq T$, the random vector $\left(Y_{t_1},Y_{t_2},\cdots,Y_{t_n}\right)$ has a joint normal distribution. From \eqref{eq:distOUP}, we know that, for any $n\geq 1$, and arbitrary $k \in \{1,\cdots,n\}$, $Y_{t_k}$ is normally distributed. From \cite[p.\ 259]{Bauer.1996}, we know that, to verify the joint normal distribution of $\left(Y_{t_1},Y_{t_2},\cdots,Y_{t_n}\right)$, it is sufficient to show that any linear combination of $Y_{t_1},Y_{t_2},\cdots,Y_{t_n}$ is normally distributed. This holds true due to the linearity of the It\^{o} integral.
	
	Furthermore, the drift coefficient $b(t,x) = \kappa \left(\mu(t) - x\right)$, and the diffusion coefficient $\sigma(t,x)\equiv \sigma$ of \eqref{eq:OUP} satisfy the assumptions for existence and uniqueness of a strong solution in \cite[Thm.\ 3.1]{Gard.1988}. Hence, \cite[Thm.\ 3.9]{Gard.1988} is applicable, which states that the corresponding SDE, in our case \eqref{eq:OUP}, is a Markov process (see \cite[Def. 4.6]{Korn.2010}) on the interval $[0,T]$.
\end{proof}

We are now in the Gauss-Markov setting of \cite{DiNardo.2001}. To prepare the numerical computation of the first passage time density, we need to calculate some basic characteristics of our OUP.
We recall, that $Y_t$ is normally distributed with mean
\begin{align*}
	m_{OUP}(t) &= y_{t_0}e^{ - \kappa (t-t_0)}  + \kappa \int\limits_{t_0}^t {e^{ - \kappa \left( {t - s} \right)} \mu \left( s \right)ds},
\end{align*}
and variance given by
\begin{align*}
	v_{OUP}(t) &= \frac{\sigma ^2}{2\kappa}\left(1 - e^{ - 2\kappa \left( {t - t_0} \right)}\right),
\end{align*}
both being $C^1\left([0,T]\right)$-functions.
The probability density function of the OUP starting at time $t_0$ in $y_{t_0}$ coincides with a normal density with mean $m_{OUP}(t)$, and variance $v_{OUP}(t)$.

We proceed with the covariance function. Note that the covariance is determined by the stochastic integral term $I_t = \sigma \int_{t_0}^{t} e^{-\kappa\left(t-s\right)}dW_s$ in \eqref{eq:expSolOUP} and the deterministic part can be neglected for its calculation. Moreover, note that $\E\left[I_t\right]=0$.
\begin{align}
	\cov(Y_s,Y_t) &= \sigma^2\E\left[\int_{t_0}^{s} e^{-\kappa\left(s-u\right)}dW_u \int_{t_0}^{t}e^{-\kappa\left(t-v\right)}dW_v\right] \notag\\
	&= \sigma^2 e^{-\kappa\left(s+t\right)} \E\left[\int_{t_0}^{s} e^{\kappa u}dW_u \int_{t_0}^{t} e^{\kappa v} dW_v\right] \notag\\
	&= \sigma^2 e^{-\kappa\left(s+t\right)} \E\left[\int_{t_0}^{t} e^{\kappa u} \ind_{[t_0,s]}(u) dW_u \int_{t_0}^{t} e^{\kappa v} dW_v\right] \notag\\
	&= \sigma^2 e^{-\kappa\left(s+t\right)} \E\left[\int_{t_0}^{s} e^{2\kappa u} du \right] \notag\\
	&= \sigma^2 \frac{e^{-\kappa\left(t-s\right)} - e^{-\kappa\left(s + t -2t_0\right)}}{2\kappa} \label{eq:covOUP}
\end{align}
The covariance function \eqref{eq:covOUP} can be decomposed in
\begin{align*}
	\cov\left(Y_s,Y_t\right) &= \underbrace{ \vphantom{\frac{\sigma^2}{2\kappa}} e^{\kappa s} \left(1 - e^{-\kappa\left(2s - 2t_0\right)}\right)}_{h_1(s)}  \cdot \underbrace{\frac{\sigma^2}{2\kappa}e^{-\kappa t}}_{h_2(t)}.
\end{align*}
Also, the functions $h1$, and $h_2$ are elements of $C^1\left([0,T]\right)$, and their derivatives are
\begin{align*}
	h_1'(t) &= \kappa e^{\kappa t} + \kappa e^{-\kappa \left(t-2t_0\right)}, \ \text{and} \\
	h_2'(t) &= -\frac{1}{2}\sigma^2e^{-\kappa t}.
\end{align*}

\bigskip
\textbf{Numerical calculation of the first passage time density}

In \cite{DiNardo.2001}, it is shown in a first step that the first passage time density for a $C^1$-barrier satisfies a non-singular Volterra second-kind integral equation. In a second step, this equation is iteratively solved by a repeated Simpson's rule yielding an approximation to the desired first passage time density in discretized form.

Below, we state Theorem 3.1.\ of \cite{DiNardo.2001} adapted to our OUP \eqref{eq:OUP} and our notation.
\begin{thm}
	Let $S(t)$ be a $C^1\left([0,T]\right)$-function. Then, the first passage time density $g\left(S(t),t|y_{t_0},t_0\right) = \frac{\partial}{\partial_t} P\left(\TT_{y_0}<t\right)$ solves the non-singualr second-kind Volterra integral equation given by
	\begin{align}
		g\left(S(t),t|y_{t_0},t_0\right) &= -2\Psi\left(S(t),t|y_{t_0},t_0\right) + 2\int_{t_0}^{t} g\left(S(t)s,s|y_{t_0},t_0\right) \Psi\left(S(t),t|S(s),s\right) ds, \quad y_{t_0}<S(t_0). \label{eq:VolterraIntEq}
	\end{align}
	Thereby, the function $\Psi$ is defined via
	\begin{align*}
		\Psi\left(S(t),t|y,s\right) =& \left(\frac{S'(t)-m_{OUP}'(t)}{2} - \frac{S(t) - m_{OUP}(t)}{2} \frac{h_1'(t)h_2(s) - h_2'(t)h_1(s)}{h_1(t)h_2(s) - h_2(t)h_1(s)} \right.\\
		&\left.- \frac{y - m_{OUP}(t)}{2} \frac{h_2'(t)h_1(t) - h_2(t)h_1'(t)}{h_1(t)h_2(s) - h_2(t)h_1(s)}\right)\cdot p_{y,s}(S(t),t).
	\end{align*}
\end{thm}
We adopt the notational short cuts introduced in \cite[p.\ 466]{DiNardo.2001}:
\begin{align*}
	g(t) &:= g\left(S(t),t|y_{t_0},t_0\right), \quad t,t_0 \in [0,T], t_0<t \\
	\Psi(t) &:= \Psi\left(S(t),t|y_{t_0},t_0\right), \quad t,t_0 \in [0,T], t_0<t \\
	\Psi\left(t|s\right) &:= \Psi\left(S(t),t|S(s),s\right) \quad t,s \in [0,T], t_0<s \leq t.
\end{align*}
We introduce a grid $t_0<t_1<...<T_N$, where $t_k = t_0 + k\cdot \Delta t, k \in \{1,\cdots,N\}$.
The iterative procedure based on the repeated Simpson's rule to obtain an approximation $\tg(t_k)$ of the first passage time density $g(t_k)$ reads as follows:
\begin{align}
	\tg(t_1) &= -2\Psi(t_1), \notag\\
	\tg(t_k) &= -2\Psi(t_k) + 2\Delta t\sum_{j=1}^{k-1} w_{k,j} \tg(t_j)\Psi(t_k|t_j), \quad k=2,3,\cdots,N \label{eq:algoFPTD}
\end{align}
where the weights are specified by
\begin{align*}
	w_{2n,2j-1} &= \frac{4}{3}, \quad j=1,2,\cdots,n; n=1,2,\cdots,\frac{N}{2}, \\
	w_{2n,2j} &= \frac{2}{3}, \quad j=1,2,\cdots,n-1; n=2,3,\cdots,\frac{N}{2}, \\
	w_{2n+1,2j-1} &= \frac{4}{3}, \quad j=1,2,\cdots,n-1; n=2,3,\cdots,\frac{N}{2}-1, \\
	w_{2n+1,2j} &= \frac{2}{3}, \quad j=1,2,\cdots,n-2; n=3,4\cdots,\frac{N}{2}-1, \\
	w_{2n+1,2(n-1)} &= \frac{17}{24}, \quad n=2,3,\cdots,\frac{N}{2}-1, \\
	w_{2n+1,2n-1} &= w_{2n+1,2n} = \frac{9}{8}, \quad n=1,2,\cdots,\frac{N}{2}-1.
\end{align*}
The iterative procedure has been proven to converge in \cite{DiNardo.2001}.
\begin{thm}
	We shall be given the above discretization $t_0<t_1<...<t_N$, where $t_k = t_0 + k\cdot \Delta t, k \in \{1,\cdots,N\}$, where $\Delta t$ is the discretization step. The first passage time density obtained by the iterative procedure \eqref{eq:algoFPTD} converges to the true first passage time density as the step size tends to zero, i.e.
	\begin{align*}
		\lim\limits_{\Delta t \rightarrow 0} |g(t_k)-\tg(t_k)| = 0 \quad \text{for all} \ k \in \{1,\cdots,N\} 
	\end{align*}
\end{thm}

To use the iteratively approximated first passage time density to obtain the risk level corresponding to $S(t)$, we apply Algorithm \ref{alg:JCC}.
\begin{algorithm}
	\caption{Algorithm to calculate the risk level corresponding to $S(t)$}
	\label{alg:JCC}
	\begin{algorithmic}[1]
		\REQUIRE OUP characteristics: mean $m_{OUP}(t)$, variance $v_{OUP}(t)$, covariance decomposition in terms of $h_1(t),h_2(t)$, and its derivatives $h_1'(t),h_2'(t)$; boundary $S(t)$; discretization step $\Delta t$, final active time $T_{CC}$ of CC; risk level $\risk$
		\ENSURE First passage time density $g$, and cumulative distribution function $\cdfJCC$
		\STATE Define parameters of OUP.
		\STATE Choose time discretization $\Delta t$.
		\STATE Calculate $m_{OUP}(t),h_1(t),h_2(t),S(t),m_{OUP}'(t),h_1'(t),h_2'(t),S'(t)$ for discretized time interval $[0,T_{CC}]$ with step size $\Delta t$.
		\STATE Calculate approximation of first passage time density via \eqref{eq:algoFPTD}.
		\STATE Calculate the corresponding discrete values of the cumulative distribution function $\cdfJCC$.
		\IF{$\cdfJCC(end)<= \risk$}
		\STATE JCC fulfilled.
		\ELSE
		\STATE JCC violated.
		\ENDIF
		\RETURN $g$ and $\cdfJCC$.
	\end{algorithmic}
\end{algorithm}
It enables to integrate the JCC \eqref{eq:JCC} in our optimization procedure to solve the SOC problem \eqref{eq:SOC-CC}. In our case the optimal inflow control will take the role of the $C^1$-boundary $S(t)$. As this control is a result of the optimization procedure, the calculation of the first passage time density needs to be repeated in every optimization iteration. It is therefore worthwhile to mention that the above introduced iterative procedure is well-suited for computational efficiency. This is because the algorithm only requires the characteristics of the OUP in terms of its initial data $(t_0,y_0)$, its mean $m_{OUP}(t)$, its variance $v_{OUP}(t)$, its covariance decomposition in terms of the functions $h_1$ and $h_2$ as well as a prespecified boundary $S(t)$ and a chosen discretization step $\Delta t$. No Monte Carlo (MC) methods, and no high-dimension integral computations are involved and no particular software packages are necessary (see \cite{DiNardo.2001}).

\subsubsection*{Validation of first passage time simulation}
For a validation of Algorithm \ref{alg:JCC}, we consider the following test case setting:
$t_0=0, \kappa=\nicefrac{1}{3600}, \sigma=0.003, \mu(t) = 0.7 + 0.3 \cdot \sin(\nicefrac{1}{7200}\ \pi t), y_0=0.8, S(t_i)=m_{OUP}(t_k) + 0.2 + 0.25\cdot \nicefrac{t_k}{T}$, where $t_k, k\in \{1,\cdots,N_{\Delta t}\}$ are the discretization points, and $N_{\Delta t}+1$ denotes the number of discretization points for the interval $[0,T]$ that corresponds to the chosen time step $\Delta t$. We denote by $M$ the number of MC repetitions used.

In Table \ref{tab:MCvalidationJCCalgorithm}, we show the risk of hitting the boundary $S(t)$ based on a MC simulation (risk-MC) and based on the evaluation of the cumulative distribution function $G$ from Algorithm \ref{alg:JCC} (risk-fptd), and calculate the differences (diff). We observe that our algorithm \ref{alg:JCC} gives rather precise results already for large step sizes. This is beneficial when using it within the optimization of large gas networks in Subsection \ref{subsec:OCgasPower}.
The MC-risk values are obtained using the plain MC method. Note the \textit{rare event character} of undersupply. This is even more pronounced in real world settings where most likely the chosen risk tolerance is $5\%$ or lower instead of values around $15\%$ in our test case. Therefore, a \textit{rare event simulation technique} as e.g.\ in \cite{Wadman.2016} in the context of power flow reliability, where the probability of an outage is very small, might lead to more accurate risk estimates even for larger step sizes.

However, even with the plain MC method, we observe that the MC-risk and the risk-fptd values approach up to a precision in the range of $10^3$ indicating the correct functioning of the algorithm.
\begin{table}[tbp]
	\setcounter{mpfootnote}{\value{footnote}}
	\renewcommand{\thempfootnote}{\arabic{mpfootnote}}
	\centering
	\caption{Comparison of risk-MC and risk values based on Algorithm \ref{alg:JCC}}
	\begin{tabular}{|rrrrrr|}
		\hline
		\multicolumn{1}{|l}{\textbf{T}} & \multicolumn{1}{l}{\textbf{$\Delta t$}} & \multicolumn{1}{l}{\textbf{M}} & \multicolumn{1}{l}{\textbf{risk-MC}} & \multicolumn{1}{l}{\textbf{risk-fptd}} & \multicolumn{1}{l|}{\textbf{diff}} \\
		\hline
		\multicolumn{1}{|l}{$4\cdot3600$} & 480   & 1000  & 0.086 & 0.1481 & 0.0621 \\
		&       & 10000 & 0.0883 & 0.1481 & 0.0598 \\
		&       & 100000 & 0.0901 & 0.1481 & 0.058 \\
		\hline
		& 60    & 1000  & 0.124 & 0.1479 & 0.0239 \\
		&       & 10000 & 0.1214 & 0.1479 & 0.0265 \\
		&       & 100000 & 0.1201 & 0.1479 & 0.0278 \\
		&       & 1000000 & 0.1198 & 0.1479 & 0.0281 \\
		\hline
		& 6     & 1000  & 0.131 & 0.1479 & 0.0169 \\
		&       & 10000 & 0.1399 & 0.1479 & 0.008 \\
		&       & 100000 & 0.1384 & 0.1479 & 0.0095 \\
		\hline
		& 1     & 1000  & 0.149 & 0.1479 & -0.0011 \\
		&       & 10000 & 0.1409 & 0.1479 & 0.007 \\
		\hline
		& 0.1   & 1000  & 0.156 & 0.1479\footnotemark[1] & -0.0081 \\
		&       & 10000 & 0.1435 & 0.1479\footnotemark[1] & 0.0044 \\
		&       & 100000 & 0.1463 & 0.1479\footnotemark[1] & 0.0016 \\
		&       & 500000 & 0.1461 & 0.1479\footnotemark[1] & 0.0018 \\
		\hline
		& 0.01  & 10000 & 0.1472 & 0.1479\footnotemark[1] & 0.0007 \\
		\hline
	\end{tabular}%
	\label{tab:MCvalidationJCCalgorithm}%
	\footnotetext[1]{Result based on step size $\Delta t=1$ for capacity reasons}
\end{table}%

\subsection{Reformulation of the cost function} \label{subsec:reformOF}
We benefit from an analytical treatment of the cost function that has been set up in \cite{Lux.2018} for the tracking type costs $C2$ of \eqref{eq:generalOF}. We use the slightly more general formulation of it allowing for arbitrary $t_0$ instead of only $t_0=0$:
\begin{align}
\E\left[\left(Y_t-S_{u(t)}\right)^2|Y_{t_0}=y_{t_0}\right] =& y_{t_0}^2e^{-2\kappa (t-t_0)} + 2y_0e^{-\kappa (t-t_0)} \int_{t_0}^{t}e^{-\kappa(t-s)}\kappa \mu(s) ds + \left(\kappa\int_{t_0}^{t}e^{-\kappa(t-s)} \mu(s) ds\right)^2 \notag\\
&+ \frac{\sigma^2}{2\kappa}\left(1-e^{-2 \kappa (t-t_0)}\right) - 2y(t) \cdot \left( y_0 e^{-\kappa (t-t_0)} + \kappa \int_{t_0}^{t} e^{-\kappa (t-s)}\mu(s) ds \right) + y(t)^2. \label{eq:reformC2}
\end{align}
It remains to extend the approach for the undersupply cost component $C3$ and the excess supply revenue component $R$. To this end, we make use of the truncated normal distribution. The following definition is taken from \cite[p.\ 81]{Johnson.1970} and adapted to our notation.
\begin{defin}
	Let $X$ be a random variable on a probability space $\left(\Omega, \AAA, P\right)$. We say that $X$ follows a \textit{doubly truncated normal distribution} with \textit{lower and upper truncation points} $a$ and $b$ respectively ($X \sim \NN_a^b(\mTGauss,\vTGauss)$) if its probability density function is given by
	\begin{align*}
		\densTGauss(x) = \left(\frac{1}{\sdevTGauss}\densGauss\left(\frac{x-\mTGauss}{\sdevTGauss}\right)\left(\Phi\left(\frac{b-\mTGauss}{\sdevTGauss}\right)-\Phi\left(\frac{a-\mTGauss}{\sdevTGauss}\right)\right)^{-1}\right)\ind_{[a,b]}(x),
	\end{align*}
	where $\densGauss$ is the density, and $\Phi$ the cumulative distribution function of a standard normally distributed random variable. We denote by $\mTGauss$ and $\vTGauss$ the mean and the variance of the non-truncated normal distribution.
	
	We call the distribution \textit{singly truncated from above} respectively \textit{from below} if $a$ is replaced by $-\infty$  respectively $b$ is replaced by $\infty$.
\end{defin}

\begin{prop}[{cf.\ \cite[p.\ 81]{Johnson.1970}}]
	Let $X \sim \NN_a^b(\mTGauss,\vTGauss)$. Then, the expected value of $X$ reads as
	\begin{align}
		\E\left[X\right] = \mTGauss + \frac{\varphi\left(\frac{a-\mTGauss}{\sdevTGauss}\right) - \varphi\left(\frac{b-\mTGauss}{\sdevTGauss}\right)}{\Phi\left(\frac{b-\mTGauss}{\sdevTGauss}\right) - \Phi\left(\frac{a-\mTGauss}{\sdevTGauss}\right)} \sdevTGauss.	\label{eq:expTruncGauss}
	\end{align}
\end{prop}
We use equation \eqref{eq:expTruncGauss} denoting the expectation of a truncated normally distributed random variable to reformulate the expectations in $C3$ and $R$ in \eqref{eq:generalOF}. For $C3$, we obtain
\begin{align}
	\E\left[\left(\Supp - Y_t\right)_{-}|Y_{t_0}=y_{t_0}\right] &= \int_{\Supp}^{\infty} \left(\Supp-y\right)\rho_{Y_t}(y) dy \notag\\
	&= \Supp \cdot \left(1-\Phi_{Y_t}(\Supp)\right) - \int_{\Supp}^{\infty} y\rho_{Y_t}(y)dy \notag\\
	&= \Supp \cdot \left(1-\Phi_{Y_t}(\Supp)\right) - \int_{\R} y\rho_{Y_t,\Supp,\infty}(y)dy \cdot P\left(Y_t > \Supp\right) \notag\\
	&= \Supp \cdot \left(1-\Phi_{Y_t}(\Supp)\right) - \left(m_{OUP}(t) + \frac{\rho_{Y_t}(\Supp)}{1-\Phi_{Y_t}(\Supp)}\sqrt{V_{OUP}(t)}\right) \cdot \left(1-\Phi_{Y_t}(\Supp)\right), \label{eq:reformC3}
\end{align}
where $\rho_{Y_t}$ denotes the density, and $\Phi_{Y_t}$ the cumulative distribution function of $Y_t$, and $\rho_{Y_t,\Supp,\infty}$ denotes the density of the singly truncated random variable from below by $\Supp$.
The last equation results from formula \eqref{eq:expTruncGauss} with $a=\Supp$, and $b=\infty$.

In a similar manner, by setting $a=-\infty$, and $b=\Supp$ in \eqref{eq:expTruncGauss}, for the expectation in $R$, we have
\begin{align}
\E\left[\left(\Supp - Y_t\right)_{+}|Y_{t_0}=y_{t_0}\right]
&= \Supp \Phi_{Y_t}(\Supp) - \left(m_{OUP}(t) - \frac{\rho_{Y_t}(\Supp)}{\Phi_{Y_t}(\Supp)}\sqrt{V_{OUP}(t)}\right)  \Phi_{Y_t}(\Supp). \label{eq:reformR}
\end{align}
With equations \eqref{eq:reformC2}, \eqref{eq:reformC3}, and \eqref{eq:reformR}, we have a completely deterministic reformulation of our cost function \eqref{eq:generalOF} at hand and denote it by $OF_{\text{detReform}}(Y_t, t_0,y_{t_0} ,\Supp)$.
This reformulation of the cost function together with the reformulation of the CC as state constraint \eqref{eq:stateConstr} allows us to drop the OUP \eqref{eq:OUP} as constraint in the optimization problem \eqref{eq:SOC-CC}.

We are left with a completely deterministic PDE-constrained optimization problem, i.e. \eqref{eq:SOC-CC} without the OUP \eqref{eq:OUP} as constraint and with \eqref{eq:generalOF} replaced by $OF_{\text{detReform}}(Y_t, t_0,y_{t_0} ,\Supp)$, and \eqref{eq:CC} replaced by \eqref{eq:stateConstr}. Hence, we are free to apply a suitable method from deterministic PDE-constrained optimization of our choice to solve the problem. Here, we use a first-discretize-then-optimize approach, and numerically calculate the discrete optimal solution based on discrete adjoints.
All together, we come up with the deterministic reformulation of the stochastic optimal control (SOC) problem given by
\begin{align}
\min_{u \in \UU_{ad}} \int_{t^{\ast}}^{T} OF_{\text{detReform}}(Y_t; t_0;y_{t_0} ;\Supp) dt \quad \text{subject to} \ \eqref{eq:IBVP}, \eqref{eq:stateConstr},\ \text{and/or} \ \eqref{eq:reformJCCasFPT}. \label{eq:reformSOC-CC}
\end{align}
Note that the JCC \eqref{eq:reformJCCasFPT} is still a probabilistic constraint, which can, however, be handled in a deterministic way by applying Algorithm \ref{alg:JCC}.

\section{Numerical results} \label{sec:numRes}
We apply deterministic discrete adjoint calculus (see e.g.\ \cite{Gottlich.2010, Goatin.2016, Kolb.2012}) to solve the deterministically reformulated SOC problem \eqref{eq:reformSOC-CC}.
Numerically, this is implemented in ANACONDA, a modularized simulation and optimization environment for hyperbolic balance laws on networks that allows the inclusion of state constraints, which has been developed in \cite{Kolb.2011}. Note that the inclusion of state constraint is important for the inclusion of the SCC \eqref{eq:stateConstr}. Further note that ANACONDA allows for regularization terms in the objective function to tackle finite horizon effects or oversteering for coarse discretizations. Both effects can be reduced by penalizing control variations, which has been activated within the computations in Sections~\ref{subsec:OClinAdvSource} and~\ref{subsec:OCtelEq} with penalty parameter $10^{-5}$. To be more precise about the applied schemes, we use the IBOX scheme \eqref{eq:IBOXscheme} to discretize the energy system \eqref{eq:IBVP}, and the IPOPT solver \cite{Wachter.2006} as well as the \textit{DONLP2} method \cite{Spellucci.1998b,Spellucci.1998a} within ANACONDA for the optimization procedure.

\bigskip
We start this section with a validation of our numerical routine in Subsection \ref{subsec:OClinAdvSource} for a special case of \eqref{eq:IBVP}, for which we are able to derive an analytical solution of the SOC problem in case of an SCC. In Subsections \ref{subsec:OCtelEq} and \ref{subsec:OCgasPower}, we apply our numerical routine in case of two different parameter settings. The choices in the energy system \eqref{eq:IBVP}, the uncertain demand \eqref{eq:OUP}, and the CCs \eqref{eq:CC} stated in Table \ref{tab:paraTrans-GtP} correspond to the settings in Sections \ref{subsec:OCtelEq}, and \ref{subsec:OCgasPower}.
\begin{table}[h!]
		\centering
		\begin{adjustbox}{max width=\textwidth}
		\begin{tabular}{lcrrrr}
			\hline 
			\textbf{Parameter setting (PS)} & & \textit{Tele} \ref{subsec:OCtelEq} & \textit{GtP\textunderscore s} \ref{subsec:OCgasPower} & \textit{GtP\textunderscore l} \ref{subsec:OCgasPower} \\ 
			\hline 
			Time horizon & $T$ &  4 & 12 & 14400 \\
			Initial demand & $y_0$ & $1$ & $0.8$ & $0.8$ \\
			Mean demand level & $\mu(t)$ & $1 + \sin(\pi t)$ & $y_0+0.3\sin(\frac{1}{7200} \pi t)$ & $0.7+0.3\sin(\frac{1}{7200} \pi t)$ \\
			Speed of mean reversion & $\kappa$ & 3 & 3 & \nicefrac{1}{3600} \\
			Intensity of demand fluctuations & $\sigma$ & 0.2 & 0.2 & 0.003 \\
			\hline
			Dimension & $n$ & 2 & 2 & 2 \\
			\# of prescribed BCs left & $m_l$ & 1 & 1 & 1 \\
			\# of prescribed BCs right & $m_r$ & 1 & 1 & 1 \\
			Functional relation of left BC & $\Gae$ & $\rho_2^e$ & $p(\rho_{1}^{e})$ & $p(\rho_{1}^{e})$ \\
			Functional relation of right BC & $\Gbe$ & $\rho_1^e$ & $\rho_2^e$ & $\rho_2^e$ \\
			Flux function & $\f(\rho^{e})$ & $\begin{pmatrix}
			-C^{-1}\rho_{2}^{e} \\
			L^{-1}\rho_{1}^{e}\\
			\end{pmatrix}$& $\begin{pmatrix}
			\rho_{2}^{e}\\
			p(\rho_{1}^{e}) + \nicefrac{\left(\rho_{2}^{e}\right)^{2}}{\rho_{1}^{e}}
			\end{pmatrix}$ & $\begin{pmatrix}
			\rho_{2}^{e}\\
			p(\rho_{1}^{e}) + \nicefrac{\left(\rho_{2}^{e}\right)^{2}}{\rho_{1}^{e}}
			\end{pmatrix}$\\
			\hline
			CC active & $I_{CC}$ & $[1.5,4]$ & $[6,12]$ & $[0,4]$ \\
			Risk level & $\risk$ & $5\%$ & $5\%$ & $5\%$ \\
			\hline
		\end{tabular}
	\end{adjustbox}
	\caption{Specification of energy systems}
	\label{tab:paraTrans-GtP}
\end{table}

\subsection{Validation via scalar linear advection with source term} \label{subsec:OClinAdvSource}
For a numerical validation of our optimization routine, we consider
the special case of the linear advection on one edge with velocity $\lambda \in \R$, and nonzero, linear source term $\s\rho$. As we only consider the dynamics on one edge, we omit the dependence on the edge $e$. The latter equation results from \eqref{eq:IBVP} by setting the dimension $n=1$, the number of prescribed left, and right BCs to $m_l=1$, and $m_r=0$, $\aae=0$ and by choosing the functional relation of the left BC as $\Ga = \id$, where $\id$ denotes the identity function, and by setting the flux function to $\f(\rho) = \lambda \rho$. It reads as
\begin{align}
\partial_t \rho + \lambda \partial_x \rho &= \s \rho \notag\\
\rho(x,0) &= \rho_0(x), \notag\\
\rho(0,t) &= u(t), \quad x \in [0,b], \ t \in [0,T]. \label{eq:linAdvWithSource-IBVP}
\end{align}
This IBVP has the advantage that we are able to derive an analytical solution to the corresponding SOC problem \eqref{eq:reformSOC-CC} with SCC \eqref{eq:stateConstr} and can compare our numerical implementation against it.

\subsubsection*{Analytical solution of the linear advection with source term}
To derive the analytical solution of the SOC problem, we first need the analytical solution of the IBVP \eqref{eq:linAdvWithSource-IBVP}. Therefore, we divide the time interval $[0,T]$ in two subintervals $I_0 = [0,\nicefrac{x}{\lambda}]$, and $I_b = (\nicefrac{x}{\lambda},T]$ depending on the position $x \in [0,b]$. The solution on $I_0$ is determined by the initial condition $\rho_0(x)$, and the solution on $I_b$ is determined by the inflow control $u(t)$ acting as a left BC.

\begin{figure}
	\centering
	\includegraphics[width=0.7\linewidth]{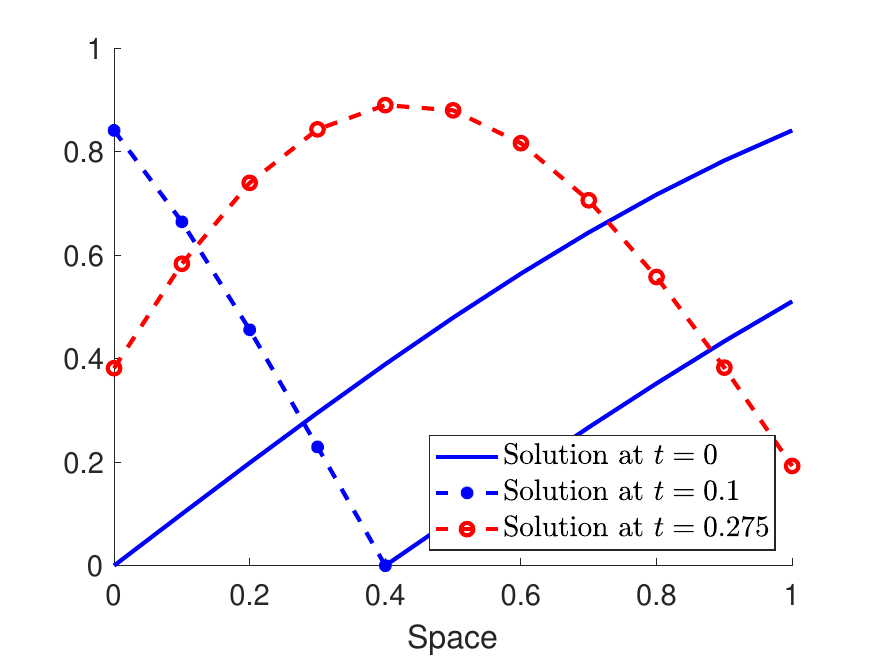}
	\caption{Regions of influence of IC ($I_0$) and inflow control ($I_b$) in solution \eqref{eq:expSolLinAdvWithSource}}
	\label{fig:influenceICandBClinAdv_source}
\end{figure}

For the solution of \eqref{eq:linAdvWithSource-IBVP} on $I_b$, we can change the role of $x$, and $t$ as we deal with an empty system at the beginning. Then, we can use the explicit solution of the IVP for the linear advection equation with source term (see \cite{Teuber.2017}).
All together, the solution of \eqref{eq:linAdvWithSource-IBVP} reads as
\begin{align}
\rho(x,t) = \left\{\begin{array}{ll}
e^{\s t} \rho_0(x-\lambda t) & t \in I_0\\
e^{\s \nicefrac{x}{\lambda}} u(t-\nicefrac{x}{\lambda}) & t \in I_b \label{eq:expSolLinAdvWithSource}
\end{array}\right.
\end{align}
The solution of \eqref{eq:linAdvWithSource-IBVP} is illustrated graphically in Figure \ref{fig:influenceICandBClinAdv_source} based on the choices $\lambda=4$, $\Delta x=0.1$, $\Delta t = \nicefrac{\Delta x}{\lambda}$ (exact CFL-condition), $T=1$, $s=-1$, $\rho_0(x)=\sin(x)$, and $u(t) = \sin(10\cdot t)$.
Note that the solution at $t=0$ is the prescribed initial condition. For $x \in [0,0.4)$, $t=0.1$ is an element of $I_b$, i.e.~the inflow control determines the solution, and for $x \in [0.4,1]$, $t=0.1$ belongs to $I_0$, which means that, thereon, the solution results from the shifted initial condition.
The solution at $t=0.275$ depends only on the inflow control.

For the particular case of \eqref{eq:linAdvWithSource-IBVP}, we are able to derive an analytical expression for the optimal control for the particular cost function $OF_{\text{detReform}}$ with $\wdet=\wunder=\wex=0$, and $\wtrack=1$.
\begin{thm} \label{thm:anaOClinAdvSource}
	Let the supply dynamics \eqref{eq:IBVP} be given by \eqref{eq:linAdvWithSource-IBVP}, and the uncertain demand by the OUP \eqref{eq:OUP} with existing second moment. Furthermore, set $\wdet=\wunder=\wex=0$, and $\wtrack=1$ in \eqref{eq:generalOF}.
	Then solving \eqref{eq:SOC-CC} with \eqref{eq:CC} being an active SCC \eqref{eq:SCC} on the whole interval $[0,T]$ yields the optimal control
	\begin{align}
	u^\ast(t) = m_{OUP}(t) + v_{OUP}(t)\Phi^{-1}\left(1-\risk\right) \ \forall t \in [0,T-\invLa] \label{eq:anaSol_linAdvSourceOneEdgeSCC}
	\end{align}
\end{thm}
\begin{proof}
	From \cite{Lux.2018}, we know that the optimal supply without imposing a CC in the linear advection setting \eqref{eq:linAdvWithSource-IBVP} with $s=0$ for CM1 is given by the conditional expectation $\Supp = \E[Y_{t+\invLa}|Y_0]$. As the source term, does not influence the characteristics of the linear advection (see equation \eqref{eq:expSolLinAdvWithSource}), we again have a constant time delay between inflow and outflow. That enables us to still consider the solution pointwise in time. Due to the analytical solution \eqref{eq:expSolLinAdvWithSource}, the optimal supply is always reachable as long as $\Supp \in [u_{\text{min}},u_{\text{max}}]$, where $u_{\text{min}}$, and $u_{\text{max}}$ are lower and upper bounds on the control ($u_{\text{min}}=0$, and $u_{\text{max}}=+\infty$ in \eqref{eq:SOC-CC}). Since we have for the chance constraint level
	\begin{align*}
	m_{OUP}(t+\invLa) + v_{OUP}(t+\invLa)\Phi^{-1}\left(1-\risk\right) \geq m_{OUP}(t+\invLa) = \E[Y_{t+\invLa}|Y_0],
	\end{align*}
	and the objective function is quadratic with vertex in $(t,\E[Y_{t+\invLa}|Y_0])$, the optimal control is obtained at the left boundary of the feasible region at each point in time.
\end{proof}

\subsubsection*{Numerical validation}
We now validate our numerical procedure by considering the particular SOC problem given by
\begin{align}
\min_{u(t), t \in [0,T], u \in L^2} \int_{t^{\ast}}^{T} OF_{\text{detReform}}(Y_t; t_0;y_{t_0} ;\Supp) dt \quad \text{subject to} \ \eqref{eq:SCC}, \ \text{and} \ \eqref{eq:linAdvWithSource-IBVP}, \label{eq:reformSOC_linAdvSource_SCC}
\end{align}
where we set $\wdet=\wunder=\wex=0$, and $\wtrack=1$ in \eqref{eq:generalOF}. We consider the interval $[0,1]$. For this setting, we derived the analytical optimal control in Theorem \ref{thm:anaOClinAdvSource} such that we can compare our numerical solution against the analytical one.

\begin{figure}
	\centering
	\setlength{\fwidth}{0.8\textwidth}
	\input{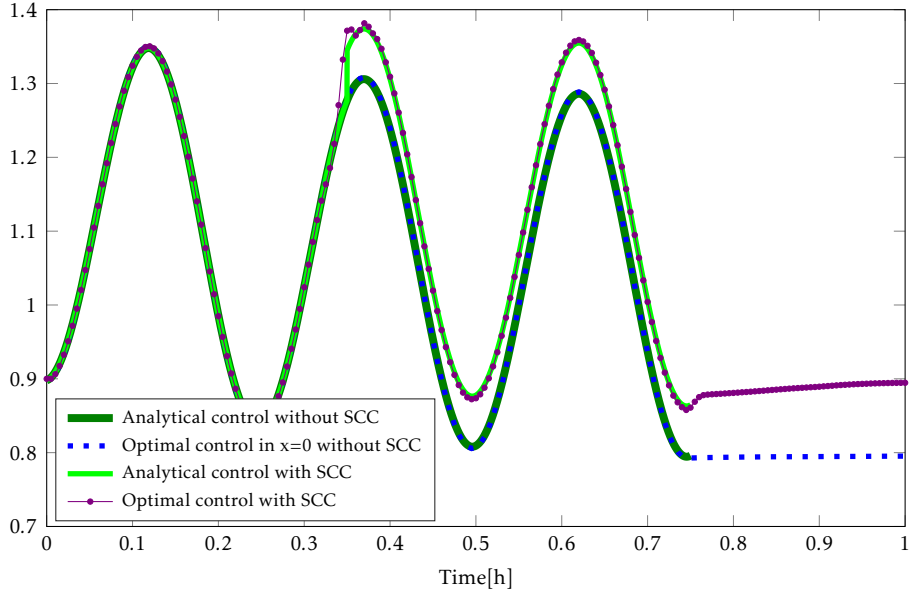}
	\caption{Numerical versus analytical solution of \eqref{eq:reformSOC_linAdvSource_SCC}}
	\label{fig:numVsAnalytical_linAdvSourceOneEdge}
\end{figure}

We activate the SCC within $I_{CC}=[0.6,1]$. The transport velocity is $\lambda=4$ and we consider a rate for the source term of $s=-0.1$.
In the numerical implementation, we have included a nonnegativity constraint $u(t)\geq 0$ on the inflow control. Note however that this constraint does not affect our test case (see Figure \ref{fig:numVsAnalytical_linAdvSourceOneEdge}). Therefore, a direct comparison of our numerical solution against the analytical one is possible.

In our validation, we used the parameters $y_0=1$, $\mu(t)=1 + 2\sin(8\pi t)$, $\kappa=3$ and $\sigma=0.1$ for the OUP.
In Figure \ref{fig:numVsAnalytical_linAdvSourceOneEdge}, we observe that the numerical optimal control is close to the analytical control \eqref{eq:anaSol_linAdvSourceOneEdgeSCC}. This finding suggests the correct functioning of our routine. Note that the small spike at the activation time of the SCC in the numerical optimal control is inherited from the jump in the analytical solution.

\subsection{Linear system of hyperbolic balance laws: Telegrapher's equation } \label{subsec:OCtelEq}
Having seen that our numerical routine provides very good approximations to the analytical solution in the case of the linear advection with source term, in this section, we now consider the Telegrapher's equation on a network, and apply our routine thereon. The corresponding energy system is obtained by specifying the quantities in the IBVP \eqref{eq:IBVP} according to the values in Table \ref{tab:paraTrans-GtP} according to the parameter setting \textit{Tele} \ref{subsec:OCtelEq}. In the following assertions, we stick to the conventional notation of the Telegrapher's equation and denote $(\rhoe_1,\rhoe_2)^{T} = (\Ue,\Ie)^{T}$, where $\Ue$ is the voltage, and $\Ie$ the current on edge $e$. We deal with a system of linear hyperbolic ba\-lance laws with linear source terms, where we assume the same constant parameters $R,L,C,G$ for resistance, inductance, capacitance, and conductance on each edge.
We endow the system with initial conditions (ICs) $\Ue_0$, and $\Ie_0$. On the bounded domain $X = [a,b]$, where $a = \min\{\aae \ |\ e \in \EE\}$, and $b = \max \{\be \ |\ e \in \EE\}$, we prescribe boundary conditions (BCs)  $v_{\text{ext}}(t)$ as an externally given function to prescribe the voltage at the end of the line, and $u(t)$ as the voltage control at the beginning of the line, taking the role of the inflow control into the energy system.
\begin{subequations} \label{eq:teleIBVP}
	\begin{align}
	\partial_t \Ue + C^{-1}\partial_x \Ie &= -GC^{-1}\Ue \notag\\
	\partial_t \Ie + L^{-1} \partial_x \Ue &= -RL^{-1}\Ie, \label{eq:telegraphers} \\
	\Ue(x,0) &= \Ue_0(x), \quad \Ie(x,0) = \Ie_0(x), \label{eq:ICstelegraphers} \\
	\Ue(a,t) &= u(t), \quad \Ue(b,t) = v_{\text{ext}}(t), \label{eq:BCstelegraphers}
	\end{align}
\end{subequations}

\bigskip

Before we present our numerical results for the Telegrapher's equation, we first address the question of well-posedness of the system on one edge, and in a second step the well-posedness of the system in the network case as well as the existence of an optimal control.

Note that the system \eqref{eq:telegraphers} is diagonalizable. Hence, we can rewrite it in characteristic variables denoted by $\xi=\left(\xi^{+},\xi^{-}\right)^T$ (see \cite{Schillen.2016}). In case of lossless transmission ($R=G=0$), we trace the system back to a decoupled system of two classical linear advection equations by splitting the dynamics into left- and right-traveling waves with characteristic speeds $\lambda^{-}$ and $\lambda^{+}$. In the lossless case for one edge normed to a length of $1$, the IBVP \eqref{eq:teleIBVP} reads as
\begin{align*}
\partial_t \xi^{+} + \lambda^{+} \partial_x \xi^{+} &= 0	\\
\partial_t \xi^{-} + \lambda^{-} \partial_x \xi^{-} &= 0 \\
\xi^{+}(x,0) = 0, \quad \xi^{-} &= 0 \\
\xi^{+}(0,t) = g^{+}(t), \quad \xi^{-}(1,t) &= g^{-}(t),
\end{align*}
where $\lambda^{\overset{+}{-}} = \overset{+}{-}\left(\sqrt{LC}\right)^{-1}$.
Note that voltage $U$, and current $I$ can be expressed in terms of the characteristic variables as $U(x,t) = \sqrt{\frac{L}{C}}\left(\xi^{+} - \xi^{-}\right)$, and $I(x,t) = \xi^{+} + \xi^{-}$.

\bigskip
As we deal with an empty system at the beginning, we can change the role of $x$, and $t$, and obtain for $\lambda^{+}>0$, and $\lambda^{-}<0$
\begin{align*}
\nicefrac{1}{\lambda^{+}} \partial_t \xi^{+} + \partial_x \xi^{+} &= 0	\\
\nicefrac{1}{\lambda^{-}} \partial_t \xi^{-} + \partial_x \xi^{-} &= 0 \\
\xi^{+}(0,t) = g^{+}(\xi^{+}(1,t)), \quad \xi^{-}(1,t) &= g^{-}(\xi^{-}(0,t)).
\end{align*}
whose solution is
\begin{align*}
\xi^{+}(x,t) &= \xi^{+}(0,t-\invLa{+}x) = g^{+}(t-\nicefrac{1}{\lambda^{+}}x) \\
\xi^{-}(x,t) &= \xi^{-}(1,t-\invLa{-}x) = g^{-}(t-\nicefrac{1}{\lambda^{-}}x) \\
\end{align*}
By interchanging the role of $x$ and $t$ again (empty system at the beginning), the existence of a unique weak solution of \eqref{eq:teleIBVP} with nonzero source term is ensured as a special case of the IVP for the system of hyperbolic balance laws with dissipative source term studied in \cite{Christoforou.2006}. For further details, we refer the reader to \cite{Teuber.2017}.

\bigskip

As we consider the dynamics given by the Telgrapher's equation \eqref{eq:telegraphers} on each edge of our network, we need to impose suitable coupling conditions at the inner nodes $v \in \VV_{int} = \VV \backslash \left(\VV_{in} \cup \VV_{out}\right)$.
Therefore, we define the set of all edges connected to node $v \in \VV$ as
$$\EEv = \{\ev_1,\cdots,\ev_{k_{\text{in}}}\} \cup \{\ev_{k_{\text{in}} + 1},\cdots,\ev_{k_{\text{in}} + k_{\text{out}}}\}.$$
We define the coupling function as
$$\couplingv: \R^{\n \times \left(|\delta^{+}(v)|+|\delta^{-}(v)|\right)} \rightarrow \R^{m_c}, c\left(\qe, e \in \delta^{-}(v) \cup \delta^{+}(v)\right) = \begin{pmatrix}
U^{\ev_1}(b^{\ev_1},t) - U^{\ev_{2}}(b^{\ev_{2}},t)\\
\vdots \\
U^{\ev_1}(b^{\ev_1},t) - U^{\ev_{k_{\text{in}}}}(b^{\ev_{k_{\text{in}}}},t)\\
U^{\ev_1}(b^{\ev_1},t) - U^{\ev_{k_{\text{in}} + 1}}(a^{\ev_{k_{\text{in}} + 1}},t)\\
U^{\ev_1}(b^{\ev_1},t) - U^{\ev_{k_{\text{in}} + k_{\text{out}}}}(a^{\ev_{k_{\text{in}}} + k_{\text{out}}},t)\\
\sum_{i \in \delta^{-}(v)} I(b^{i},t) - \sum_{j \in \delta^{+}(v)} I(a^{j},t)
\end{pmatrix}$$
for each vertex $v \in \VV_{int}$. Thereby, $m_c$ denotes the number of coupling conditions. Now, we simply state the coupling conditions as $c_v = 0$ for all $v \in \VV_{int}$.
\begin{rem}
	The coupling conditions with respect to $\Uev, \ev \in \EEv,$ state the equality of the voltages at the corresponding boundaries of all edges with common node $v$. The coupling condition for $\Ie$ at node $v$ models the conservation of the flow of current.
\end{rem}

We are interested in the existence of an optimal inflow control $u(t)$ for the system \eqref{eq:teleIBVP}.
Similar problems have been tackled in \cite{Herty.2009}. There, they derive the well-posedness of a system of nonlinear hyperbolic balance laws on a network with edges represented by the interval $[0,\infty)$, and the existence of an optimal control minimizing a given cost functional under certain conditions in the context of gas networks and open canals.

As our flow function $\fe(\xi) = \Lambda\xi$ is linear with $\lambda^{-} < 0 < \lambda^{+}$, it satisfies the assumption \textbf{(F)} in \cite{Herty.2009}. Moreover, our source term $\s(\xi) = \B$ is a constant function implying its Lipschitz property and its bounded total variation. Hence, condition \textbf{(G)} from \cite{Herty.2009} is also fulfilled. This gives us the well-posedness of our IBVP \eqref{eq:teleIBVP} on a network on the positive half-line (see \cite[Theorem 2.3.]{Herty.2009}).
%

\bigskip

We now turn our attention to the full SOC problem \eqref{eq:reformSOC-CC} in consideration.
An existence result of an optimal control in a gas network setting in the presence of state constraints has been derived in \cite{Wollner2017}, where the state constraints enter the optimization problem via a barrier method. Note that the SCC \eqref{eq:SCC} for our SOC problem can be incorporated as a state constraint. The existence of a minimizer in the presence of a JCC \eqref{eq:JCC} is more involved as a reformulation as state constraint is no longer possible.

Therefore, we focus on the numerical analysis of the SOC problem \eqref{eq:SOC-CC} with particular focus on the influence of different types of chance constraints.
\subsubsection*{Telegrapher's equation with single chance constraint} \label{subsubsec:CC-OC-tele-net}
We start with a single chance constraint \eqref{eq:SCC} in the Telegrapher's setting on the network depicted in Figure \ref{fig:networkTop} with parameters specified in Table \ref{tab:paraTrans-GtP} (\textit{Tele} \ref{subsec:OCtelEq}), and consider the cost function \eqref{eq:generalOF} with $\wdet=\wunder=\wex=0$, and $\wtrack=1$. We consider the network topology from Figure \ref{fig:networkTop}. The initial conditions are set to $\Ue(x,0)=1$ on all edges and $\Ie(x,0) = 1$ for $e_1$ and $e_7$, $\Ie(x,0) = 0.5$ for $e_2$, $e_3$, $e_5$ and $e_6$ as well as $\Ie(x,0) = 0$ for $e_4$. We consider $[\aae,\be]=[0,1]$ for edges $e \in \{e_1,\dots,e_6\}$ and $[\aae,\be] = [0,2]$ for $e_7$.  We set the model parameters in \eqref{eq:teleIBVP} to $R=0.01$ (resistance), $L=0.5$ (inductance), $C=0.125$ (capacitance), and $G=0.01$ (conductance).

\begin{figure}
	\centering
	\setlength{\fwidth}{0.7\textwidth}
	\input{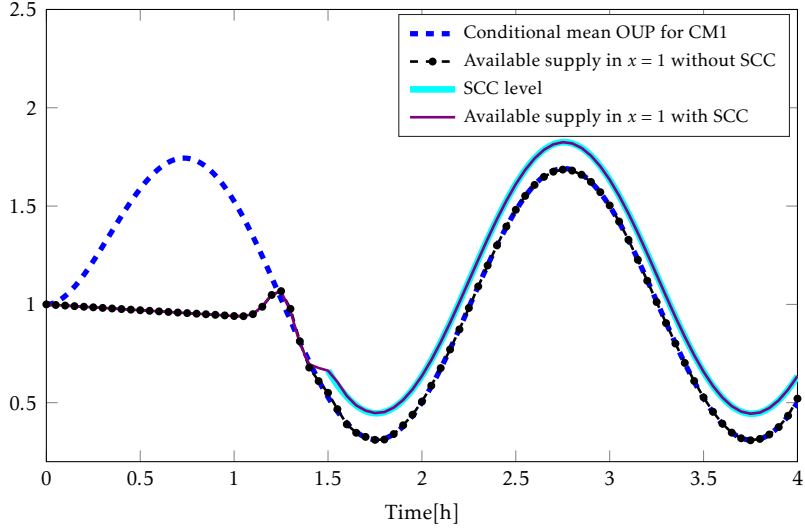}
	\caption{Optimal supply for Telegrapher's equation on network of Figure \ref{fig:networkTop}}
	\label{fig:gammtelegraphersnetwork}
\end{figure}

Our numerical results in Figure \ref{fig:gammtelegraphersnetwork} show that the optimally available current $I$ (purple solid line) at node $v_d$ matches the expected value of the OUP (blue dashed line) till $t = 1.5$, and follows the course of the optimally available current without SCC (black dotted line). For $t \in I_{CC}$, the optimal available current matches the given CC-Level (turquoise solid line). This goes along with our intuitive understanding: the optimal available current matches the $95\%$- confidence level. This is because a higher supply results in an even higher tracking-type error in the cost function, and a lower supply violates the SCC at the risk level of $5\%$. This numerical finding is in line with the theoretical result in Theorem \ref{thm:anaOClinAdvSource}.

\subsection{Nonlinear system of hyperbolic balance laws: gas-to-power system} \label{subsec:OCgasPower}
The energy transition phase comes with a lot of challenges in its realization. Aiming for a high percentage of energy provided from renewable energy sources, one has to somehow cope with the large volatility in the energy generation.
One possibility to react to fluctuations to still ensure a stable demand satisfaction are gas turbines. The possibility of gas-to-power, i.e. the withdrawal of gas from the gas system and its transformation to power, seems to be a promising approach. A gas turbine can be booted to full performance within several minutes. Thus, a gas turbine might be appropriate to overcome short-term bottlenecks in energy generation.
Therefore, we finally focus on the modeling of gas-to-power in the context of uncertain demands. We particularly pay attention how the withdrawal of gas affects the behavior of the gas network.

One major mathematical challenge that the gas network entails is that its governing equations are no longer a linear system of hyperbolic balance laws.
The nonlinear gas transport can be described by the so-called \textit{isentropic Euler equations}. They can be obtained from equation \eqref{eq:IBVP} by putting in the parameters of \textit{GtP} \ref{subsec:OCgasPower} of Table \ref{tab:paraTrans-GtP}:
\begin{subequations}
\begin{align}
		\begin{pmatrix} \rho^{e} \\ q^{e} \end{pmatrix}_t &+ \begin{pmatrix} q^{e} \\ p(\rho^{e}) + \frac{(q^{e})^2}{\rho^{e}} \end{pmatrix}_x = \begin{pmatrix} 0 \\ g(\rho^{e},q^{e}) \end{pmatrix} \label{eq:isentropicEulerEq}\\
		\rho^{e}(x,0) &= \rho_0^{e}(x), \quad q^{e}(x,0) = q_0^{e}(x).		\label{eq:isentropicEulerICBC}
\end{align}
\end{subequations}
Thereby, $\rhoe_1$ is denoted by $\rho^{e}$ and describes the density, $\rhoe_2$ is identified with $q$ and gives the flow, $g$ a given source term chosen as in \cite{Fokken.2019}, and $\pe$ is the pressure on edge $e$, which is specified by the pressure law $\pe(\rho^{e})=d^2\cdot\rho^{\beta}$ for all edges $e$ (with $\beta=1$ and $d=340 \frac{\text{m}}{\text{s}}$ in the examples below). Note that the exponent $\beta=1$ means that we deal with the isothermal Euler equations, a special case of the isentropic Euler equations. However, the analysis is not limited to this choice and results for various pressure laws can for example be found in \cite{Fokken.2018}.

\begin{figure}
	\centering

	\begin{tikzpicture}[scale=0.5, every node/.append style={transform shape}]
	\node[draw, circle, fill=RoyalBlue] (vd) at (0,0) {\textbf{$v_d$}};  
	
	\node (powerStation) at (vd.north)[xshift =-7.3cm ,yshift=2cm] {\includegraphics[width=0.25\linewidth]{myPowerStation.png}};
	\node (demand)  at (powerStation.east)[xshift =8.8cm ,yshift=-0.5cm]{\includegraphics[width=0.4\linewidth]{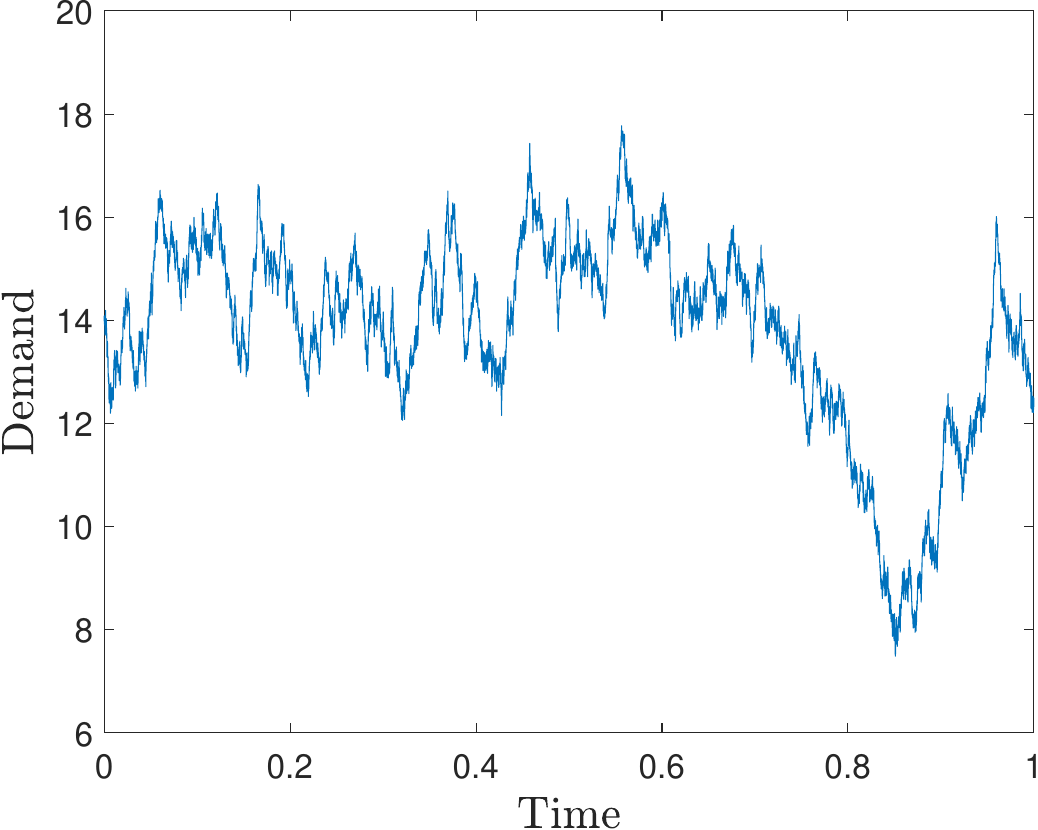}};
	\node at (demand)[xshift =-0.3cm ,yshift=-0.7cm]{\large \textcolor{RoyalBlue}{$Y_t$}};	
	\node (extSupp) at (powerStation.south)[xshift =0.2cm ,yshift=-0.2cm] {\Large \textbf{External supply:}};
	\node at (extSupp)[xshift =0cm ,yshift=-0.7cm]{\Large $\tilde{Y}_t=Y_t-S_{u(t)}$};
	\node (gasSupp) at (vd.south)[xshift =-2.8cm ,yshift=-0.5cm] {\Large \textbf{Supply through gas}};
	\node at (gasSupp.south)[xshift =-0.6cm ,yshift=-0.5cm] {\Large \textbf{$q=\epsilon(S_{u(t)})$}};
	\node at (gasSupp.south)[xshift =0.4cm ,yshift=-1.3cm] {\Large \textbf{$=a_0 + a_1 S_{u(t)} + a_2 S_{u(t)}^2$}};
	\node[draw=RoyalBlue, rectangle, thick, minimum width=16.7cm,
	minimum height = 7.8cm] at (vd)[xshift =-1.3cm ,yshift=1cm] {};  
	\node[draw, circle] (vin) at (vd.south)[xshift =-8cm ,yshift=-5.2cm] {$v_{in}$};
	\node[draw, rectangle] (v0) at (vin.east)[xshift =2cm ,yshift=0cm] {$v_0$};  
	\node[draw, circle] (v1) at (vin.east)[xshift =5cm ,yshift=0cm] {$v_1$};
	\node[draw, circle] (v4) at (vin.east)[xshift =10cm ,yshift=0cm] {$v_4$};  
	\node[draw, circle] (v5) at (vin.east)[xshift =13cm ,yshift=0cm] {$v_5$};  
	\node[draw, circle] (v2) at (vd.south)[xshift =0cm ,yshift=-3.2cm] {$v_2$};  
	\node[draw, circle] (v3) at (v2)[xshift =0cm ,yshift=-4cm] {$v_3$};
	
	\draw[-latex, thick] (vin) to node[above,yshift=0.1cm] {\Large $e_0$} (v0);
	\draw[-latex, line width=1.5pt] (v0) to node[above,yshift=0.1cm] {\Large $e_1$} (v1);
	\draw[-latex, thick] (v1) to node[above,yshift=0.1cm] {\Large $e_2$} (v2);
	\draw[-latex, thick] (v1) to node[above,yshift=0.1cm] {\Large $e_3$} (v3);
	
	\draw[-latex, thick] (v2) to node[above,yshift=0.1cm] {\Large $e_4$} (v4);
	\draw[-latex, thick] (v3) to node[above,yshift=0.1cm] {\Large $e_5$} (v4);
	\draw[-latex, thick] (v4) to node[above,yshift=0.1cm] {\Large $e_6$} (v5);
	
	\draw (v0.east)[xshift =1cm ,yshift=-0.6cm] circle (8pt);
	\node[draw, cross out] (compressor)  at (v0.east)[xshift =1cm ,yshift=-0.6cm]{};
	\node  at (compressor.north)[xshift =0cm ,yshift=-0.8cm]{\Large compressor};
	\node[draw=red!20, rectangle, thick, minimum width=16.7cm,
	minimum height = 5cm] at (vd)[xshift =-1.3cm ,yshift=-5.6cm] {};

	\draw[-latex, line width=2pt, dashed] (v2) edge (vd);
	\draw[-latex, line width=2pt, dashed] (powerStation) edge (vd);
	\end{tikzpicture}

	\caption{Coupled gas-to-power system}
	\label{fig:gasToPower}
\end{figure}

To ensure the well-posedness of the system, we also need coupling conditions at the nodes. As in \cite{Fokken.2019}, we use pressure equality and mass conservation at all nodes at all times (Kirchhoff-type-coupling), i.e. for all $v \in \VV$ and for all $t \in [0,T]$, we have
\begin{align*}
	p_{in}^v(t) &= p_{out}^v(t), \\
	\sum_{e \in \delta^{-}(v)} \qe(t) &= \sum_{\tilde{e} \in \delta^{+}(v)} q^{\tilde{e}}(t).
\end{align*}
To couple the gas network to the power system, we adapt the deterministic coupled gas-to-power-system described in \cite{Fokken.2019}, and extend it by an uncertain power demand given by the OUP \eqref{eq:OUP}. The adapted setting is sketched in Figure \ref{fig:gasToPower}. 
In our case, the power system is shrinked to only one node $v_d$, where the aggregated uncertain demand is realized. This is because our focus is on matching this demand $Y_t$ best possible by the power $S_{u(t)}$ provided by the gas-to-power turbine, whereas the distribution within the power system is considered as a separate task. 

The power demand is attained preferably by the gas-to-power conversion amount $S_{u(t)}$. The missing power $\tilde{Y}_t$ needs to be covered by an external power source.

The gas consumption to generate the power is described by a quadratic function
\begin{align*}
\epsilon(\Supp) = a_0 + a_1\Supp + a_2(\Supp)^2.
\end{align*}
As the amount of gas withdrawn from the network is our control, we have
\begin{align*}
	u(t) = \epsilon(\Supp).
\end{align*}


Pressure bounds play an important role in gas networks. In our case, we add a lower pressure bound at node $v_5$, which is
\begin{align}
p^{v_5}(t) \geq 43 [bar] \label{eq:pressureBoundv5}
\end{align}
in the first scenario below, appearing as an additional algebraic constraint in the SOC \eqref{eq:reformSOC-CC}.

\bigskip
There is a pressure drop in the gas network caused by the gas withdrawal for the conversion to power. It can be compensated by a compressor station. The modeling of the compressor station is taken from \cite{Fokken.2019}: the compressor is modeled as edge $e_1$ with time-independent in- and outgoing pressure $p^{v_0}$ and $p^{v_1}$, and flux values $q^{v_0}$, and $q^{v_1}$ through the nodes $v_0$ and $v_1$. We assume that the compressor is run via an external power source only linked to the scenario in Figure \ref{fig:gasToPower} via the cost component $C1$ in \eqref{eq:generalOF} accounting for operator costs of the compressor. This entails the flux equality $q^{v_1}=q^{v_0}$. Note that the operating costs increase if the ratio $\nicefrac{p^{v_1}}{p^{v_0}}$ increases. The compressor can be controlled by the gas network operator. Therefore, a second control variable, i.e.\ the control of the compressor station $\compr$, is added to the SOC problem \eqref{eq:reformSOC-CC} in terms of 
\begin{align*}
	\compr &= p^{v_1}(t) - p^{v_0}(t).
\end{align*}
It appears in the operating costs $C1^{\text{compr}}$ in the cost component $C1$ of $OF_{\text{detReform}}$.
Moreover, costs occur for the gas consumption to satisfy the power demand. These costs complete the deterministic costs $C1$ in \eqref{eq:generalOF} as
\begin{align*}
C1 = C1^{\text{compr}} + 0.0001\cdot \int_{t_0}^{T} u(t) dt.
\end{align*}
For our purpose, it is important to note the deterministic nature of the compressor costs $C1^{\text{compr}}$ and the costs for the gas consumption within the objective function \eqref{eq:generalOF}.

\bigskip
Due to the active element in terms of the compressor station, we need to adapt two of the above coupling conditions at $v_1$ and $v_2$ as
\begin{subequations}
	\begin{align}
	p_{out}^{v_1}(t) &= p_{in}^{v_1}(t) + \compr, \ \text{and} \label{eq:contrCompr}\\
	q_{out}^{v_2}(t) &= q_{in}^{v_2}(t) - u(t). \label{eq:withdrawGas}
	\end{align}
\end{subequations}
Equation \eqref{eq:contrCompr} describes the gas coupling accounting for the possible pressure increase at $v_1$. The mass conservation at $v_2$ except for the gas withdrawn for the gas-to-power conversion is modeled by equation \eqref{eq:withdrawGas}.

\bigskip

It might well be the case that the power demand cannot completely be served by the gas network due to pressure bounds in combination with the maximal performance of the compressor station, or other physical or economical reasons. In this case, the external power supplier comes into play. Hence, the difference $Y_t-\Supp$ is covered by external power supply.

\subsubsection*{Gas-to-power system with single chance constraint}
For the numerical investigation of the impact of an SCC on the optimal amount of gas withdrawn from the network, we consider the parameter setting GtP\textunderscore s in Table \ref{tab:paraTrans-GtP}. In the cost function \eqref{eq:generalOF} of the SOC \eqref{eq:reformSOC-CC}, we consider a pure tracking type functional by setting $\wdet=\wunder=\wex=0$, and $\wtrack=1$.

Our numerics refer to the setting schematically represented in \ref{fig:gasToPower}. The gas network is a subgrid of the GasLib-40 network, which approximates a segment of the low-calorific gas network located in the Rhine-Main-Ruhr area in Germany \cite{Schmidt.2017}. An extension of the gas network by a compressor station in a purely deterministic setting has already been analyzed in \cite{Fokken.2019}. Here, we present results including the compressor station in the presence of an uncertain demand stream as well as an SCC. To the best of our knowledge, this has not been considered before. We choose a discretization of $dt=900[\text{s}]$ and $dx \approx 1[\text{km}]$, slightly adapted to the individual length of the edges. The specifications of the intervals $[\aae,\be]$ for the edges are due to the network specifications in \cite{Schmidt.2017} and can be found in \cite[Table 1]{Fokken.2019} for the subgrid considered here.
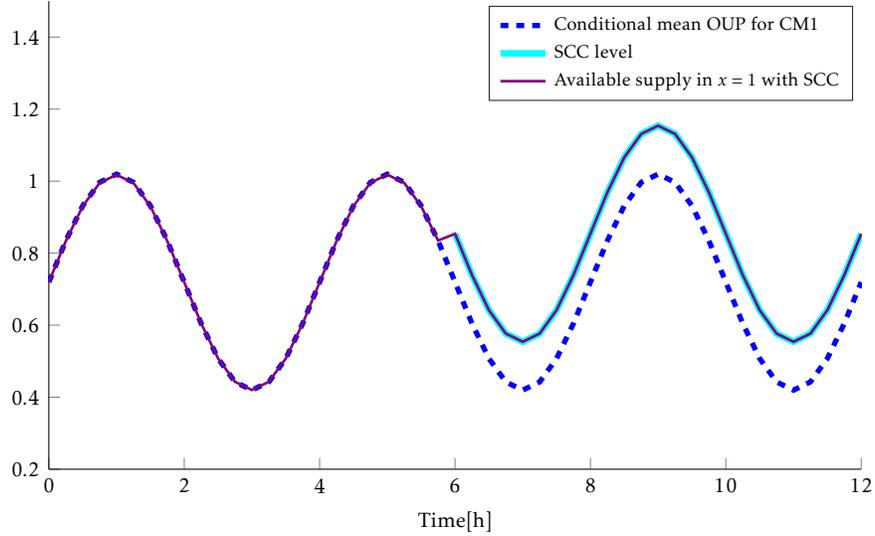
\begin{figure}[tbp]
	\centering
	\setlength{\fwidth}{0.8\textwidth}
%
%
\definecolor{mycolor1}{rgb}{0.00000,1.00000,1.00000}%
\begin{tikzpicture}

\begin{axis}[%
width=0.9\fwidth,
height=0.519\fwidth,
at={(0\fwidth,0\fwidth)},
scale only axis,
point meta min=0,
point meta max=1,
xmin=0,
xmax=12,
xlabel style={font=\color{white!15!black}},
xlabel={Time[h]},
ymin=0.2,
ymax=1.5,
axis background/.style={fill=white},
axis x line*=bottom,
axis y line*=left,
label style={font=\footnotesize},
tick label style={font=\footnotesize},
legend style={at={(0.99,0.99)}, anchor=north east, legend cell align=left, align=left, draw=white!15!black, font = \scriptsize},
]

\addplot [color=blue, dashed, line width=2.0pt]
  table[row sep=crcr]{%
0	0.719547\\
0.25	0.834311715432461\\
0.5	0.931648176515441\\
0.75	0.9966941561759\\
1	1.0195469936538\\
1.25	0.99672755160463\\
1.5	0.931709883221611\\
1.75	0.834392339129426\\
2.25	0.604782284567539\\
2.5	0.507445823484559\\
2.75	0.4423998438241\\
3	0.419547006346196\\
3.25	0.442366448395369\\
3.5	0.50738411677839\\
3.75	0.604701660870575\\
4.25	0.834311715432461\\
4.5	0.931648176515441\\
4.75	0.9966941561759\\
5	1.0195469936538\\
5.25	0.99672755160463\\
5.5	0.931709883221611\\
5.75	0.834392339129426\\
6.25	0.604782284567539\\
6.5	0.507445823484559\\
6.75	0.4423998438241\\
7	0.419547006346196\\
7.25	0.442366448395369\\
7.5	0.50738411677839\\
7.75	0.604701660870575\\
8.25	0.834311715432461\\
8.5	0.931648176515441\\
8.75	0.9966941561759\\
9	1.0195469936538\\
9.25	0.996727551604632\\
9.5	0.931709883221611\\
9.75	0.834392339129424\\
10.25	0.604782284567539\\
10.5	0.507445823484559\\
10.75	0.4423998438241\\
11	0.419547006346196\\
11.25	0.442366448395369\\
11.5	0.50738411677839\\
11.75	0.604701660870575\\
12	0.719503366769622\\
};
\addlegendentry{Conditional mean OUP for CM1}

\addplot [color=mycolor1, line width=2.5pt]
table[row sep=crcr]{%
	6	0.853892369483622\\
	6.25	0.739084020820785\\
	6.5	0.641747559737803\\
	6.75	0.576701580077344\\
	7	0.55384874259944\\
	7.25	0.576668184648614\\
	7.5	0.641685853031634\\
	7.75	0.73900339712382\\
	8.25	0.968613451685705\\
	8.5	1.06594991276869\\
	8.75	1.13099589242914\\
	9	1.15384872990705\\
	9.25	1.13102928785788\\
	9.5	1.06601161947486\\
	9.75	0.96869407538267\\
	10.25	0.739084020820783\\
	10.5	0.641747559737803\\
	10.75	0.576701580077346\\
	11	0.55384874259944\\
	11.25	0.576668184648614\\
	11.5	0.641685853031634\\
	11.75	0.73900339712382\\
	12	0.853805103022868\\
};
\addlegendentry{SCC level}

\addplot [color=violet, line width=1pt]
  table[row sep=crcr]{%
0	0.719547\\
0.25	0.833320000000001\\
0.5	0.930265\\
0.75	0.995002\\
1	1.017697\\
1.25	0.994945\\
1.5	0.930294\\
1.75	0.833621000000001\\
2	0.719481999999999\\
2.25	0.604896\\
2.5	0.507736\\
2.75	0.442805\\
3	0.419979\\
3.25	0.442726\\
3.5	0.507574999999999\\
3.75	0.60464\\
4	0.719131000000001\\
4.25	0.833603\\
4.5	0.930621\\
4.75	0.995404000000001\\
5	1.018097\\
5.25	0.995257000000001\\
5.5	0.930408\\
5.75	0.835032\\
6	0.853892999999999\\
6.25	0.739084\\
6.5	0.641748\\
6.75	0.576702000000001\\
7	0.553849\\
7.25	0.576669000000001\\
7.5	0.641686\\
7.75	0.739004\\
8	0.853804999999999\\
8.25	0.968614000000001\\
8.5	1.06595\\
8.75	1.130996\\
9	1.153849\\
9.25	1.13103\\
9.5	1.066012\\
9.75	0.968693999999999\\
10	0.853892999999999\\
10.25	0.739084\\
10.5	0.641748\\
10.75	0.576702000000001\\
11	0.553849\\
11.25	0.576669000000001\\
11.5	0.641686\\
11.75	0.739004\\
12	0.853804999999999\\
};
\addlegendentry{Available supply in $x=1$ with SCC}

\end{axis}
\end{tikzpicture}%
	\caption{Optimal amount of gas-to-power conversion}
	\label{fig:mcwoptdemand}
\end{figure}
\begin{figure}[tbp]
	\centering
	\setlength{\fwidth}{0.4\textwidth}
	\subfloat[\ Pressure increase at compressor station \label{fig:woptpressureincreasecompressor}]{
%
%
\begin{tikzpicture}

\begin{axis}[%
width=0.951\fwidth,
height=0.573\fwidth,
at={(0\fwidth,0\fwidth)},
scale only axis,
xmin=0,
xmax=12,
xlabel style={font=\color{white!15!black}},
xlabel={Time[h]},
ymin=-0.0191265,
ymax=0.4016565,
ylabel style={font=\color{white!15!black}},
ylabel={Pressure increase [bar]},
axis background/.style={fill=white},
axis x line*=bottom,
axis y line*=left,
label style={font=\footnotesize},
tick label style={font=\footnotesize}
]
\addplot [color=black, line width=2.0pt, forget plot]
  table[row sep=crcr]{%
0	0\\
0.25	0.0481200000000008\\
0.5	0.0889380000000006\\
0.75	0.117967999999999\\
1	0.130796999999999\\
1.25	0.124005\\
1.5	0.0970359999999992\\
1.75	0.0557700000000008\\
2	0.0165539999999993\\
2.25	0.000605000000000189\\
2.5	0.000199999999999534\\
3.5	0.000308000000000419\\
3.75	0.000405999999999906\\
4	0.00136799999999937\\
4.25	0.00762399999999985\\
4.5	0.0187889999999999\\
4.75	0.0320459999999994\\
5	0.0441079999999996\\
5.25	0.0515270000000001\\
5.5	0.0511789999999994\\
5.75	0.0414849999999998\\
6	0.0247670000000007\\
6.25	0.0103460000000002\\
6.5	0.0110650000000003\\
6.75	0.0274800000000006\\
7	0.0573429999999995\\
7.25	0.0979700000000001\\
7.5	0.146314\\
7.75	0.198959\\
8	0.25215\\
8.25	0.301779\\
8.5	0.343356999999999\\
8.75	0.371999000000001\\
9	0.382529999999999\\
9.25	0.369837\\
9.5	0.329917999999999\\
9.75	0.262357\\
10	0.174884\\
10.25	0.0886309999999995\\
10.5	0.0313890000000008\\
10.75	0.0025429999999993\\
11	0.000159999999999272\\
11.75	0.000147000000000119\\
12	0.00028599999999912\\
};
\end{axis}
\end{tikzpicture}%
}
	\subfloat[\ Pressure at the sink \label{fig:wopt_pressureincreasecompressor}]{
%
%
\definecolor{mycolor1}{rgb}{0.85000,0.32500,0.09800}%
\begin{tikzpicture}

\begin{axis}[%
width=0.951\fwidth,
height=0.573\fwidth,
at={(0\fwidth,0\fwidth)},
scale only axis,
xmin=0,
xmax=12,
xlabel style={font=\color{white!15!black}},
xlabel={Time[h]},
ymin=42.95,
ymax=43.4,
ylabel style={font=\color{white!15!black}},
ylabel={Pressure [bar]},
axis background/.style={fill=white},
axis x line*=bottom,
axis y line*=left,
label style={font=\footnotesize},
tick label style={font=\footnotesize}
]
\addplot [color=blue, line width=1.0pt, forget plot]
  table[row sep=crcr]{%
0	43.228301\\
0.25	43.2203263985588\\
0.5	43.1974768176861\\
0.75	43.158924272992\\
1	43.1101411485816\\
1.25	43.0607364106287\\
1.5	43.0212248953003\\
1.75	43.000011512553\\
2	43.0018903225215\\
2.25	43.0282305224674\\
2.5	43.0764372912256\\
2.75	43.1399729721612\\
3	43.2101746783137\\
3.25	43.277974327351\\
3.5	43.3347872331961\\
3.75	43.3729071749661\\
4	43.3861008005728\\
4.25	43.3709659245383\\
4.5	43.3282602584988\\
4.75	43.2636903514662\\
5	43.1876172118163\\
5.25	43.1132770508377\\
5.5	43.0537983323257\\
5.75	43.0188846572927\\
6	43.0000200144731\\
6.25	43.0000980734197\\
6.5	43.0221879912126\\
6.75	43.0641278093426\\
7	43.1196314142423\\
7.25	43.1804158676562\\
7.5	43.2375847208132\\
7.75	43.2824253899284\\
8	43.3072647181597\\
8.25	43.3067906778105\\
8.5	43.2796727686498\\
8.75	43.2297672435954\\
9	43.1660349734356\\
9.25	43.10071165847\\
9.5	43.0461144108246\\
9.75	43.0112885129313\\
10	43.0000090273272\\
10.25	43.0112006616453\\
10.5	43.0411710668337\\
10.75	43.0845161586427\\
11	43.1348929216736\\
11.25	43.1843088913713\\
11.5	43.2238276741669\\
11.75	43.2446997715446\\
12	43.2395735109429\\
};
\addplot [color=mycolor1, forget plot]
  table[row sep=crcr]{%
0	43\\
12	43\\
};
\addplot [color=black, line width=2.0pt, forget plot]
  table[row sep=crcr]{%
0	43\\
12	43\\
};
\end{axis}
\end{tikzpicture}%
}
	\caption{Pressure evolution}
	\label{fig:pressureEvol}
\end{figure}
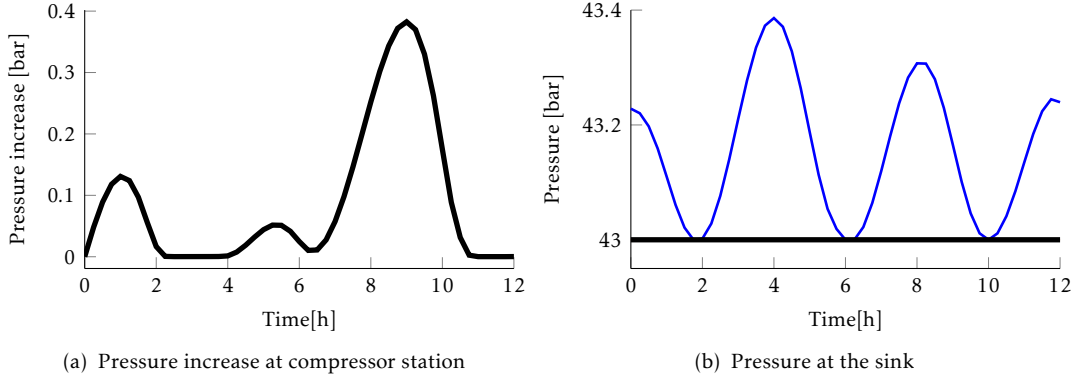

In Figure \ref{fig:mcwoptdemand}, we see that the optimal amount of gas-to-power conversion (purple solid line) follows well the course of the expected value of the OUP \eqref{eq:OUP} (blue dashed line), and from $t=6$ on matches the given SCC level (turquoise solid line) until $T=12$. This again coincides with the theoretical result for the linear advection with source term in Theorem \ref{thm:anaOClinAdvSource}.

In Figure \ref{fig:pressureEvol}, the necessity of the compressor station as an active element to keep the lower pressure bound is illustrated. There is a pressure increase at the compressor station with time-shift to ensure the lower pressure bound at $v_5$, which is set to $p^{v_5}(t) \equiv 43 [bar]$ here. This behavior has already been observed and investigated in a deterministic setting in \cite{Fokken.2019}. As we would expect, this increase is particularly pronounced while the SCC is active. The lower bound on the supply induces the need for a higher power level. This entails a higher gas withdrawal inducing a pressure drop in the gas network. The latter drop needs to be compensated by the compressor station to keep the lower pressure bound at the sink.

\subsubsection*{Gas-to-power system: comparison of single and joint chance constraint}

\begin{figure}[htbp]
	\centering
	\setlength{\fwidth}{0.6\textwidth}
	\input{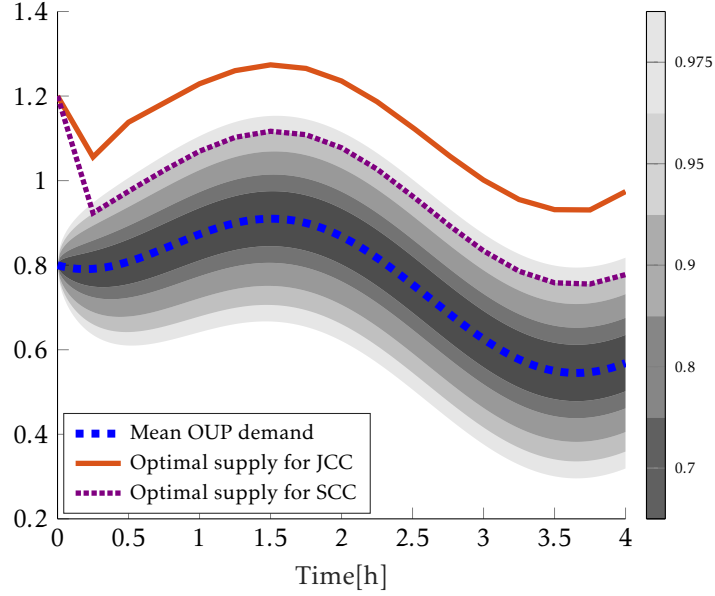}
	\caption{Optimal amount of gas to power conversion with JCC compared to SCC}
	\label{fig:optSuppJCCvsSCC_gasToPowerNetwork}
\end{figure}

In a next step, we consider the general objective function \eqref{eq:generalOF} being able to depict real costs.
To do so, we include operating costs of the compressor station (C1), costs for external power supply (C3), as well as profit of selling excess power from the gas conversion (R), and set $\wdet=\wunder=10^{-4}$, $\wtrack=0$, and $\wex=10^{-6}$. Since the cost component $C2$ cannot directly be interpreted as real costs, we set $\wtrack=0$. To handle the terms in the cost function, we use the deterministic reformulation of the additional cost components from Subsection \ref{subsec:reformOF}. Furthermore, we impose a JCC on the interval $I_{CC}=[0,4]$.
We consider a discretization of $dt=60[\text{s}]$ and $dx\approx 1[\text{km}]$, again slightly adapted to the individual length of the edges, applied to the same subgrid of the GasLib-40 network \cite{Schmidt.2017} as above. Note that this coarse discretization is possible due to the chosen numerical scheme, the IBOX scheme \eqref{eq:IBOXscheme}. The coarse time discretization is possible due to the properties of the IBOX scheme introduced in Subsection \ref{subsec:hypSupplySystem}. The control grid differs from the one in the simulation. The amount of gas-to-power conversion can be adapted every 15 minutes. We work with a lower pressure bound of $p^{v_5}(t) \equiv 42.5$.
In Figure \ref{fig:optSuppJCCvsSCC_gasToPowerNetwork}, we compare the optimal amount of gas withdrawal for the JCC with the one that would be obtained for an SCC at the same risk level of $\risk=5\%$ within the same interval $I_{CC}$.
The grey scale shows the pointwise quantile levels of the OUP and the blue dashed line indicates the mean of the OUP.
Our numerical results reveal that the optimal supply for the JCC (red solid line) evolves in the same structural manner as the SCC (purple dotted line) but at a higher level. This is in line with our intuitive understanding: the JCC acts pathwise and thus represents a more restrictive constraint. Note that a path that once violated the constraint cannot be counted as \textit{save path}, i.e.~a path below the imposed CC, at any later point in time anymore. However, this is possible for an SCC.

\begin{figure}[htbp]
	\centering
	\setlength{\fwidth}{0.4\textwidth}
	\subfloat[\ Evolution of the running maximum \label{fig:runningMaxJCC}]{\input{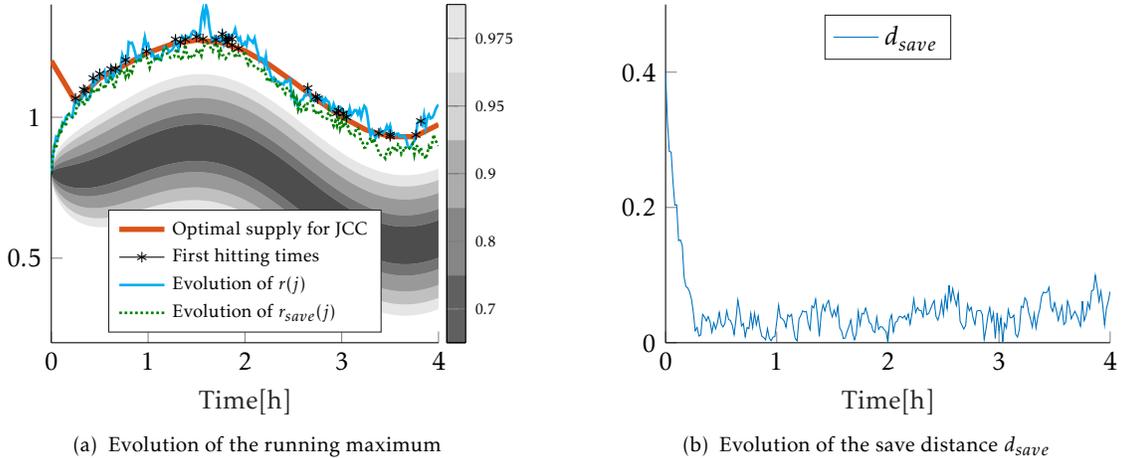}}\hfill
	\subfloat[\ Evolution of the save distance $d_{save}$ \label{fig:saveDistanceJCC}]{
%
%
\definecolor{mycolor1}{rgb}{0.00000,0.44700,0.74100}%
\begin{tikzpicture}

\begin{axis}[%
width=0.985\fwidth,
height=0.75\fwidth,
at={(0\fwidth,0\fwidth)},
scale only axis,
xmin=0,
xmax=4,
xlabel style={font=\color{white!15!black}},
xlabel={Time[h]},
ymin=0,
ymax=0.5,
axis background/.style={fill=white},
axis x line*=bottom,
axis y line*=left,
legend style={at={(0.5,0.97)}, anchor=north, legend cell align=left, align=left, draw=white!15!black}
]
\addplot [color=mycolor1]
  table[row sep=crcr]{%
0	0.4\\
0.0166666666666666	0.326054998913473\\
0.0333333333333332	0.282782809383472\\
0.0499999999999998	0.281837631147097\\
0.0666666666666664	0.237381505245147\\
0.083333333333333	0.203371322695962\\
0.0999999999999996	0.203853075170416\\
0.116666666666666	0.150963732810443\\
0.133333333333334	0.151925384247853\\
0.15	0.143917509475834\\
0.166666666666667	0.0968731049685214\\
0.183333333333334	0.0830059628895832\\
0.2	0.0809355881245724\\
0.216666666666667	0.0748087203131469\\
0.233333333333333	0.0640415191675272\\
0.25	0.0442593635172743\\
0.266666666666667	0.0178277957474524\\
0.283333333333333	0.0224218086284882\\
0.3	0.00988244516547176\\
0.316666666666666	0.0115749048769489\\
0.333333333333333	0.0420415078409144\\
0.35	0.0382013592960124\\
0.366666666666666	0.0322056564159672\\
0.383333333333334	0.0384724051991263\\
0.4	0.0312065843782721\\
0.416666666666667	0.0153793020868962\\
0.433333333333334	0.0294038252042768\\
0.45	0.0206455093134039\\
0.466666666666667	0.0295773105580919\\
0.483333333333333	0.0446477969905423\\
0.5	0.0278638066770549\\
0.516666666666667	0.0233904647414009\\
0.533333333333333	0.0281059171157665\\
0.55	0.0471486690664866\\
0.566666666666666	0.0345636777395812\\
0.583333333333333	0.0327745288054837\\
0.6	0.0261225450903693\\
0.616666666666667	0.0122995354493316\\
0.633333333333334	0.0342638341301074\\
0.65	0.0630172834848839\\
0.666666666666667	0.0473230647882579\\
0.683333333333334	0.0248434188049282\\
0.7	0.0282926561442887\\
0.716666666666667	0.00905289528545783\\
0.733333333333333	0.0224344060319934\\
0.75	0.0209403169952402\\
0.766666666666667	0.0431234906784495\\
0.783333333333333	0.0224629746475351\\
0.8	0.0439667627006033\\
0.816666666666666	0.0266511992874987\\
0.833333333333333	0.0391768781476296\\
0.85	0.0208253059057402\\
0.866666666666667	0.00921369059743515\\
0.883333333333334	0.00778517357762354\\
0.9	0.00422534380490358\\
0.916666666666667	0.0209187954215802\\
0.933333333333334	0.00665629928061051\\
0.95	0.00300191826069973\\
0.966666666666667	0.00841417545759171\\
0.983333333333333	0.00913431625247707\\
1	0.0311988796164657\\
1.01666666666667	0.0555545441519891\\
1.03333333333333	0.0370722299021073\\
1.05	0.0357038874855089\\
1.06666666666667	0.032111701020507\\
1.08333333333333	0.0230520120396429\\
1.1	0.0412705573836973\\
1.11666666666667	0.0384501852356696\\
1.13333333333333	0.0420194785516852\\
1.15	0.0211139444044033\\
1.16666666666667	0.00604241457926769\\
1.18333333333333	0.00288840819840352\\
1.2	0.0349747758659085\\
1.21666666666667	0.0404286508935465\\
1.23333333333333	0.047157850996042\\
1.25	0.0455074705566929\\
1.26666666666667	0.0463086434728455\\
1.28333333333333	0.0283103965254599\\
1.3	0.0194068903226237\\
1.31666666666667	0.0544131123727798\\
1.33333333333333	0.0516594563245532\\
1.35	0.0641382991244353\\
1.36666666666667	0.0402002707679046\\
1.38333333333333	0.00959233520878389\\
1.4	0.0427727882660625\\
1.41666666666667	0.055920018794259\\
1.43333333333333	0.0470493631519107\\
1.45	0.0574591431401741\\
1.46666666666667	0.0623658487756318\\
1.48333333333333	0.0472165018366244\\
1.5	0.0382802314448609\\
1.51666666666667	0.0481686060966062\\
1.53333333333333	0.0393416807616864\\
1.55	0.0533415407196705\\
1.56666666666667	0.0218873483838316\\
1.58333333333333	0.0172141512622828\\
1.6	0.0060304522552741\\
1.61666666666667	0.00719056594463652\\
1.63333333333333	0.0167855021362842\\
1.65	0.0240554347922917\\
1.66666666666667	0.0335679380210134\\
1.68333333333333	0.0140138031950396\\
1.7	0.00446731529381417\\
1.71666666666667	0.0138638113799709\\
1.73333333333333	0.046023835005875\\
1.75	0.0263781772205167\\
1.76666666666667	0.0434828230986026\\
1.78333333333333	0.0425444473096368\\
1.8	0.0511734215381772\\
1.81666666666667	0.0415108550751597\\
1.83333333333333	0.0296799318987118\\
1.85	0.0275604290954012\\
1.86666666666667	0.0174818140491277\\
1.88333333333333	0.0152080360672073\\
1.9	0.037650724347424\\
1.91666666666667	0.0198427191871078\\
1.93333333333333	0.023809117765361\\
1.95	0.0233126344416705\\
1.96666666666667	0.0193993146689726\\
1.98333333333333	0.016439714518496\\
2	0.000763934555190637\\
2.01666666666667	0.0327544638295265\\
2.03333333333333	0.0410585083658113\\
2.05	0.0163637312874094\\
2.06666666666667	0.0298502405552279\\
2.08333333333333	0.0316805260820701\\
2.1	0.0384808195166952\\
2.11666666666667	0.0387298304931738\\
2.13333333333333	0.0364165968321872\\
2.15	0.0398694577693073\\
2.16666666666667	0.0342288668237529\\
2.18333333333333	0.042171063457964\\
2.2	0.0274225131604302\\
2.21666666666667	0.0428457079361042\\
2.23333333333333	0.0646563279157162\\
2.25	0.04485242519727\\
2.26666666666667	0.0519073279970401\\
2.28333333333333	0.0477351784399778\\
2.3	0.0475955918139306\\
2.31666666666667	0.0384507345465233\\
2.33333333333333	0.0619021825352624\\
2.35	0.0415552421552281\\
2.36666666666667	0.0431977844395917\\
2.38333333333333	0.0495142912775055\\
2.4	0.0565816554220557\\
2.41666666666667	0.0621165549428069\\
2.43333333333333	0.0386867966600288\\
2.45	0.0334551204108964\\
2.46666666666667	0.0430493140117703\\
2.48333333333333	0.0611458459824226\\
2.5	0.0412418747793346\\
2.51666666666667	0.0692294338614667\\
2.53333333333333	0.0506251342311268\\
2.55	0.0848856638638917\\
2.56666666666667	0.061986293601505\\
2.58333333333333	0.0797293725743327\\
2.6	0.0655941152930435\\
2.61666666666667	0.0550590655895675\\
2.63333333333333	0.0726346004802787\\
2.65	0.0317969094777943\\
2.66666666666667	0.0257855702983454\\
2.68333333333333	0.0201932438960855\\
2.7	0.0254064694153602\\
2.71666666666667	0.0569053696404849\\
2.73333333333333	0.0472168225659875\\
2.75	0.0523288542667775\\
2.76666666666667	0.0364871497650663\\
2.78333333333333	0.0272301599898706\\
2.8	0.0466914364943598\\
2.81666666666667	0.0475083856154059\\
2.83333333333333	0.0476889438122212\\
2.85	0.0417666275193929\\
2.86666666666667	0.00512294628884202\\
2.88333333333333	0.0205485735516451\\
2.9	0.039840343402787\\
2.91666666666667	0.0516607938538014\\
2.93333333333333	0.0318497173304699\\
2.95	0.0438165548012739\\
2.96666666666667	0.00661931371915792\\
2.98333333333333	0.0440697287283243\\
3	0.0176071212233913\\
3.01666666666667	0.0365424812346804\\
3.03333333333333	0.00168127484794844\\
3.05	0.0283841538221461\\
3.06666666666667	0.0216423126840475\\
3.08333333333333	0.0488177862872128\\
3.1	0.0367666963550919\\
3.11666666666667	0.0371600684205786\\
3.13333333333333	0.0407495879776452\\
3.15	0.0123861976307422\\
3.16666666666667	0.00838961447155473\\
3.18333333333333	0.019687038742477\\
3.2	0.054013470951392\\
3.21666666666667	0.0519901276542978\\
3.23333333333333	0.0663165690352034\\
3.25	0.0615006207201727\\
3.26666666666667	0.0449636278746128\\
3.28333333333333	0.0353592599926715\\
3.3	0.0231408514249036\\
3.31666666666667	0.0430489243284562\\
3.33333333333333	0.0218290453067418\\
3.35	0.0390089156587194\\
3.36666666666667	0.0496874710658899\\
3.38333333333333	0.0682855059303904\\
3.4	0.0748997945795225\\
3.41666666666667	0.0738788815960261\\
3.43333333333333	0.0748402154513208\\
3.45	0.0819562749277951\\
3.46666666666667	0.0528878128421821\\
3.48333333333333	0.0463528893228453\\
3.5	0.0589571953573778\\
3.51666666666667	0.0473418222445536\\
3.53333333333333	0.0330885549729496\\
3.55	0.0443093774001939\\
3.56666666666667	0.0542805957812718\\
3.58333333333333	0.0269233516220142\\
3.6	0.0301544403196825\\
3.61666666666667	0.0613650182225696\\
3.63333333333333	0.0523380091561201\\
3.65	0.0354741776400864\\
3.66666666666667	0.0414808441390138\\
3.68333333333333	0.0465443085029147\\
3.7	0.0408481573601041\\
3.71666666666667	0.0331960661949697\\
3.73333333333333	0.0245427805787752\\
3.75	0.0613801656675328\\
3.76666666666667	0.0330506145893814\\
3.78333333333333	0.0410122535852659\\
3.8	0.0593558669184899\\
3.81666666666667	0.0506096820913173\\
3.83333333333333	0.07679575828982\\
3.85	0.0780099709741506\\
3.86666666666667	0.0994311813902993\\
3.88333333333333	0.0823547787329471\\
3.9	0.0611465980646182\\
3.91666666666667	0.0760883577478797\\
3.93333333333333	0.0552847106447736\\
3.95	0.0264056105004471\\
3.96666666666667	0.0507966974979519\\
3.98333333333333	0.0612032076682167\\
4	0.0757214045100119\\
};
\addlegendentry{$d_{save}$}

\end{axis}
\end{tikzpicture}%
}
	\caption{MC analysis of optimal supply with JCC}
	\label{fig:MCanalysisOptSuppJCC}
\end{figure}

\bigskip
We conclude our numerical investigation of the impact of chance constraints on the optimal supply with a Monte Carlo (MC) investigation of the obtained optimal supply levels for the JCC. We introduce two time-dependent sets of paths distinguishing those paths that at least once hit the prescribed supply boundary $S(t)$ until time $t$ from those that stay below this level until time $t$, i.e.
\begin{align*}
	\Omega_{hit}(t) &= \{\omega \in \Omega \ |\ \exists s \in [0,t]: Y_s(\omega)>S(t)\}, \ \text{and} \\
	\Omega_{save}(t) &= \{\omega \in \Omega \ |\ \forall s \in [0,t]: Y_s(\omega)\leq S(t)\}.
\end{align*}
Note that $\Omega_{hit} \cup \Omega_{save} = \Omega$. We denote the elements in $\Omega_{hit}$ \textit{hit paths}, and the elements in $\Omega_{save}$ \textit{save paths}. In Figure \ref{fig:runningMaxJCC}, we depict the optimal supply level with JCC as the solid blue line. The confidence intervals of the OUP are shown in grey scale. Note that they are calculated pointwise in time corresponding to an SCC. We investigate the running maximum of the save paths (green dotted line) and over all paths (light blue solid line), which we define via the evaluations of the numerical approximation $\hat{Y}_{j}$ of the OUP obtained for $M=10^3$ MC samples.
This leads to
\begin{align*}
r(j) &= \max_{ \{k \in 1,\dots,M \}} \hat{Y}_j(\omega_k), \quad j \in \{1,\dots,N_{\Delta t}\} \quad \text{and} \\
r_{save}(j) &= \max_{\{k \in 1,\dots,M \ | \ \hat{Y}_j(\omega_k) \leq S(t_j) \}} \hat{Y}_j(\omega_k), \quad \ \text{for} \ j \in \{1,\dots,N_{\Delta t}\},
\end{align*}
where $N_{\Delta t}+1$ is the number of time grid points. Note that the running maximum $r_{save}(j)$ is only well-defined if the index set is nonempty. Furthermore, we depict the instances of the first passage time of each hit paths as black asterisks.

We observe that the running maximum $r_{save}(t)$ stays close below the optimal supply $S(t)$ whereas the running maximum $r_{hit}(t)$ evolves slightly above $S(t)$. This can be interpreted in the following way: we minimize the expected quadratic deviation between our supply and the demand at the market as a major component in our cost function. At the same time, a guarantee that no undersupply occurs with a $95\%$ probability is given. Hence, it appears logical that we control our energy system in a way being close to this lower probabilistic undersupply bound to avoid costly oversupply occurrences.

\bigskip
With respect to the range where we expect the demand to evolve within (confidence intervals in grey scale), our optimal supply with JCC appears by far too large.
However, looking at the save distance $d_{save}(j) = S(t_j) - r_{save}(j)$ in Figure \ref{fig:saveDistanceJCC}, we observe primarily values below $0.1$ indicating a very limited buffer explaining the comparably large optimal supply.

\section{Conclusion}
In this work, we analyzed the optimal inflow control in hyperbolic energy systems subject to uncertain demand and particularly focused on quantifying the related uncertainty. By imposing chance constraints, we were able to limit the risk of an undersupply to a chosen probability level $\theta$. Thereby, we distinguished between SCCs and JCCs. Whereas it turned out that there is a quantile-based reformulation of the SCC enabling us to include the SCC as a state constraint into the optimization, the JCC case was more involved. We came up with a first passage time reformulation of the JCC, which we solved by using an algorithm set up in the general context of Gauss-Markov processes \cite{DiNardo.2001}.


To show that our numerical procedure can be applied to real-world phenomena, we took a real gas network from the GasLib-40 library \cite{Schmidt.2017} and calculated the optimal amount of gas withdrawn from the network to cover an uncertain power demand stream described by an OUP in the presence of a JCC.

An interesting aspect for further research is the theoretical existence of an optimal control for the above analyzed hyperbolic energy systems in the presence of a JCC.

\section*{Acknowledgement}
The authors are grateful for the support of the German Research Foundation (DFG) within the project ``\textit{Novel models and control for networked pro\-blems: from discrete event to continuous dynamics}'' (GO1920/4-1) and the BMBF within the project ``\textit{ENets}''.
Moreover, the author Kerstin Lux would like to thank Hanno Gottschalk for suggesting the consideration of joint chance constraints in this context.


\bibliographystyle{siam}
\bibliography{mysources}

\end{document}